# RADEMACHER'S THEOREM ON CONFIGURATION SPACES AND APPLICATIONS

Michael Röckner and Alexander Schied

ABSTRACT: We consider an $L^2$-Wasserstein type distance $\rho$ on the configuration space $\Gamma_X$ over a Riemannian manifold $X$, and we prove that $\rho$-Lipschitz functions are contained in a Dirichlet space associated with a measure on $\Gamma_X$ satisfying some general assumptions. These assumptions are in particular fulfilled by a large class of tempered grandcanonical Gibbs measures with respect to a superstable lower regular pair potential. As an application we prove a criterion in terms of $\rho$ for a set to be exceptional. This result immediately implies, for instance, a quasi-sure version of the spatial ergodic theorem. We also show that $\rho$ is optimal in the sense that it is the intrinsic metric of our Dirichlet form.

## 0. Introduction.

Let $\Gamma_X$ be the configuration space over a Riemannian manifold $X$. In this paper, we consider a class of probability measures on $\Gamma_X$, which in particular contains certain Ruelle type Gibbs measures and mixed Poisson measures. Using a natural 'non-flat' geometric structure of $\Gamma_X$, recently analyzed in Albeverio, Kondratiev and Röckner (1996a), (1996b), (1997a), and (1997b), one can define weak derivatives and introduce the related Sobolev spaces. Here we are interested in a more detailed description of this concept of differentiability. Similar to the case of $H$-differentiability on Wiener space, it turns out that not all values which a function $u$ takes in a small neighborhood of some $\gamma \in \Gamma_X$ are relevant for its weak gradient $\nabla^\Gamma u(\gamma)$, but only those which are located in certain 'directions'. Here, of course, the word 'direction' needs to be defined because of the absence of any linear structure on $\Gamma_X$ even if $X = \mathbb{R}^d$.

Research supported in part by Deutsche Forschungsgemeinschaft and, at MSRI, by NSF grant DMS-9701755.

Our way of making the above precise is to prove an infinite dimensional version of the celebrated theorem of Rademacher (1919) stating that Lipschitz functions on $\mathbb{R}^d$ are differentiable almost everywhere and and in the weak sense. On abstract Wiener space and its generalizations, similar results were obtained by Kusuoka (1982a), (1982b), Enchev and Stroock (1993), and Bogachev and Mayer-Wolf (1996). On configuration space, the correct Lipschitz condition is defined through an $L^2$-Wasserstein type distance function $\rho$, which, for non-compact $X$, divides $\Gamma_X$ into uncountably many disjoint 'fibers', each of the form $\{\omega \,|\, \rho(\gamma, \omega) < \infty\}$. A consequence of our Rademacher type theorem then is that only the behavior of $u$ in small $\rho$-balls around $\gamma$ matters for the value of $\nabla^\Gamma u(\gamma)$.

In a second result, we prove a partial converse to our Rademacher theorem. It also allows us to identify $\rho$ as the intrinsic metric of certain Dirichlet forms associated with our weak gradient. The resulting variational formula can also be regarded as a Kantorovich-Rubinstein type theorem for our $L^2$-Wasserstein metric. On abstract Wiener space such a converse to Rademacher's theorem was obtained in Enchev and Stroock (1993) by a different method. As a by-product to our proof we obtain that all measures satisfying our assumptions have full topological support on $\Gamma_X$, a result which might be of independent interest, in particular for the Ruelle measures mentioned above.

Another main part of this paper is devoted to applications of the above results to the potential theory on configuration space. In particular, we show that if $A \subset \Gamma_X$ has full measure, then the set of all points with positive $\rho$-distance to $A$ is exceptional. This, for instance, implies immediately a quasi-sure version of the spatial ergodic theorem on $\Gamma_{\mathbb{R}^d}$. We also give a short proof of the quasi-regularity of our Dirichlet forms.

The organization of the paper is as follows. In Section 1, we describe our setup and the Rademacher type results under some general conditions on a measure $\mu$ on $\Gamma_X$. In Section 2, we identify a class of Ruelle type Gibbs measures which satisfy these assumptions. The applications to potential and ergodic theory are presented in Section 3. Sections 4 and 5 are devoted to a detailed analysis of the topological and geometric properties of our $L^2$-Wasserstein distance $\rho$. These results serve as a preparation for the proofs of our main theorems, but some may also be of independent interest. Proofs concerning our various Sobolev spaces and Dirichlet forms are given in Section 6. Finally the proofs of our first two theorems are presented in Sections 7 and 8.

## 1. The Rademacher-type results.

Let $\Gamma_X$ denote the space of all integer-valued Radon measures on a Riemannian manifold $X$ with Riemannian inner product $g$. We will assume throughout this paper that $(X, g)$ is smooth, connected and complete. In the case where $X$ is one-dimensional we will assume that $X$ equals $\mathbb{R}$ with its Euclidean metric. Throughout this paper, we will assume that $\Gamma_X$ is endowed with the topology of vague convergence. Then $\Gamma_X$ is Polish as a vaguely closed subset of the space of all non-negative Radon measures on $X$. The set of all $\gamma \in \Gamma_X$ such that $\gamma(\{x\}) \in \{0, 1\}$ is usually called *configuration space* over $X$, but we will also call $\Gamma_X$ itself configuration space.

The space $\Gamma_X$ carries the following geometric structure which was defined in Albeverio, Kondratiev and Röckner (1996), (1997a). The "tangent space" $T_\gamma \Gamma_X$ of $\Gamma_X$ in some $\gamma$ is given as $L^2(X \to TX, \gamma)$, i.e., the space of all sections $V$ in the tangent bundle $TX$ of $X$ which are square-integrable with respect to $\gamma$:

$$\|V\|^2_{T_\gamma \Gamma_X} := \langle V, V \rangle_{T_\gamma \Gamma_X} := \int g_x(V, V) \, \gamma(dx) < \infty.$$

For ease of notation, we will also use $\|\cdot\|_\gamma$ and $\langle \cdot, \cdot \rangle_\gamma$ instead of $\|\cdot\|_{T_\gamma \Gamma_X}$ and $\langle \cdot, \cdot \rangle_{T_\gamma \Gamma_X}$. By endowing the tangent space $T_\gamma \Gamma_X$ with the inner product $\langle \cdot, \cdot \rangle_\gamma$, $\Gamma_X$ obtains a Riemannian-type structure, which is is non-trivial (i.e., varies with $\gamma$) even when the underlying space $X$ is flat.

A suitable space of "smooth test functions" on $\Gamma_X$ is the space $\mathcal{F}C_b^\infty$ which consists of all functions $u$ on $\Gamma_X$ of the form

$$(1.1) \qquad u(\gamma) = F\big(\langle f_1, \gamma \rangle, \ldots, \langle f_n, \gamma \rangle\big), \qquad \gamma \in \Gamma_X,$$

for some $n \in \mathbb{N}$, $F \in C_b^\infty(\mathbb{R}^n)$, and $f_1, \ldots, f_n \in C_0^\infty(X)$. For $u$ as in (1.1), we define its "gradient" $\nabla^\Gamma u$ as a mapping from $\Gamma_X \times X$ to $TX$, i.e., as a section of the tangent bundle $T\Gamma_X = \bigcup_\gamma T_\gamma \Gamma_X$:

$$(1.2) \qquad \nabla^\Gamma u(\gamma, x) := \sum_{i=1}^n \partial_i F\big(\langle f_1, \gamma \rangle, \ldots, \langle f_n, \gamma \rangle\big) \nabla^X f_i(x), \qquad \gamma \in \Gamma_X, \, x \in X.$$

Here $\partial_i$ means partial derivative in direction of the $i$-th coordinate, and $\nabla^X$ is the usual gradient on $X$. Alternatively, $\nabla^\Gamma u$ can be obtained using directional derivatives on $\Gamma_X$. To this end, let $V_0(X)$ denote the space of all smooth vector fields on $X$ having compact support. Then, for fixed $\gamma \in \Gamma_X$, the flow of diffeomorphisms $(\psi_t)_{t \in \mathbb{R}}$ generated by some $V \in V_0(X)$ induces a curve $\psi_t^* \gamma := \gamma \circ (\psi_t)^{-1}$, $t \in \mathbb{R}$, on $\Gamma_X$. With these notations we get that

$$\frac{d}{dt}\bigg|_{t=0} u\big(\psi_t^* \gamma\big) = \langle \nabla^\Gamma u(\gamma), V \rangle_\gamma =: \nabla^\Gamma_V u(\gamma), \qquad \gamma \in \Gamma_X, \, u \in \mathcal{F}C_b^\infty,$$

which in particular implies that (1.2) does not depend on a special representation of $u$ as in (1.1).

Let us now introduce an $L^2$-*Wasserstein type distance* $\rho$ on $\Gamma_X$ as follows.

$$\rho(\omega, \gamma) := \inf \left\{ \sqrt{\int d(x,y)^2 \, \eta(dx, dy)} \,\bigg|\, \eta \in \Gamma_{\omega, \gamma} \right\}, \qquad \omega, \gamma \in \Gamma_X,$$

where $\Gamma_{\omega,\gamma}$ denotes the set of $\eta \in \Gamma_{X \times X}$ having marginals $\omega$ and $\gamma$, and $d$ is the Riemannian distance function on $X$. Note that $\rho(\omega,\gamma)$ will be infinite if $\omega(X) \neq \gamma(X)$, because $\Gamma_{\omega,\gamma}$ will then be empty. But also if both $\omega$ and $\gamma$ are infinite configurations one will find that $\rho(\omega,\gamma) = \infty$ in general as the following example shows. Take $X = \mathbb{R}$, $\omega = \sum_{z \in \mathbb{Z}} \delta_z$, and $\gamma = \omega - \delta_0$, where $\delta_z$ denotes the Dirac measure in $z$. Obviously, convergence with respect to $\rho$ implies vague convergence.

Let us now formulate a set of assumptions we will impose in the sequel on a probability measure $\mu$ on $\Gamma_X$.

**Assumption 1.1:** We suppose that $\mu$ is a Borel probability measure on $\Gamma_X$ such that the following conditions hold.

(a) $\gamma(\{x\}) \in \{0,1\}$, for all $x \in X$ and $\mu$-a.e. $\gamma$.
(b) The mapping $\gamma \mapsto \gamma(K)$ is in $L^2(\mu)$, for each compact $K \subset X$.
(c) For any $n \in \mathbb{N}$, either $\mu(\{\gamma \,|\, \gamma(X) = n\}) > 0$ or $\mu(\{\gamma \,|\, \gamma(X) \geq n\}) > 0$ corresponding to whether $X$ is compact or non-compact respectively.
(d) For all $V \in V_0(X)$ and $t \in \mathbb{R}$, $\mu$ is quasi-invariant with respect to the flow $\psi_t^*$ of $V$, i.e., $\mu \circ (\psi_t^*)^{-1} \approx \mu$. Moreover we assume that $\mu$-a.s. $\text{ess inf}_{r \leq s \leq t} \Phi_s > 0$, for all finite $r < t$, where
$$\Phi_s := \frac{d\mu \circ (\psi_{V,s}^*)^{-1} \otimes ds}{d\mu \otimes ds}.$$

(e) $\mu$ satisfies the following integration by parts formula. If $u, v \in \mathcal{F}C_b^\infty$ and $V \in V_0(X)$, then there is an element $\nabla_V^{\Gamma*} v \in L^2(\Gamma_X, \mu)$ such that
$$\int \nabla_V^\Gamma u \, v \, d\mu = \int u \, \nabla_V^{\Gamma*} v \, d\mu.$$

Our main examples of measures satisfying Assumption 1.1 are Ruelle-type Gibbs measures corresponding to a pair potential satisfying certain assumptions (cf. Proposition 2.1 below). They will be discussed in detail in Section 2. Another example is provided by the following class of mixed Poisson measures.

**Example 1.2:** Let $m$ denote the canonical Riemannian volume element on $X$ and fix a measure $\sigma$ having a smooth and strictly positive density with respect to $m$. Then consider a measure $\mu$ on $\Gamma_X$ which is given as
$$\mu = \int_{[0,\infty)} \pi_{s \cdot \sigma} \, \lambda(ds),$$
where $\pi_{s \cdot \sigma}$ denotes the Poisson measure with intensity $s \cdot \sigma$ (for $s = 0$, $\pi_{s \cdot \sigma}$ will be the Dirac mass on the empty configuration), and $\lambda$ is a probability measure on $[0, \infty)$ such that $\int s^2 \lambda(ds) < \infty$ and $\lambda(\{0\}) < 1$. Then we claim that $\mu$ satisfies Assumption 1.1. Indeed, (a) and (c) are trivial, (b) follows from our assumptions on $\lambda$, (e) was shown in Albeverio, Kondratiev and Röckner (1997a), and quasi-invariance under diffeomorphisms is well-known with an explicit formula for the densities which then implies the positivity in (d) (see e.g. Albeverio, Kondratiev and Röckner (1997a)).

**Theorem 1.3:** *Suppose that $\mu$ satisfies Assumption 1.1 and $u \in L^2(\mu)$ is $\rho$-Lipschitz continuous. Then there exist a measurable subset $\Gamma_0$ of $\Gamma_X$ having full $\mu$-measure and a measurable section $\nabla^\Gamma u$ of $T\Gamma_X$ with the following properties.*
(i) *For all $\gamma \in \Gamma_0$, $\|\nabla^\Gamma u(\gamma)\|_\gamma \leq \mathrm{Lip}(u)$.*
(ii) *If $V \in V_0(X)$ is a vector field possessing the flow $(\psi_t)_{t \in \mathbb{R}}$, then*

$$\frac{u(\psi_t^* \gamma) - u(\gamma)}{t} \longrightarrow \langle \nabla^\Gamma u(\gamma), V \rangle_\gamma, \qquad \text{as } t \to 0,$$

*for all $\gamma \in \Gamma_0$ and in $L^2(\mu \circ (\psi_s^*)^{-1})$, for all $s \in \mathbb{R}$.*

Theorem 1.3 suggests that $\rho$-Lipschitz functions in $L^2(\mu)$ should be contained in a suitable Sobolev space. To make this idea precise, let us introduce the following Dirichlet spaces. Let $\mathbb{F}$ denote the set of all bounded measurable functions on $\Gamma_X$ for which there exists a measurable section $\nabla^\Gamma u$ of $T\Gamma_X$ such that

$$\mathcal{E}_\mu^\Gamma(u,u) := \int \langle \nabla^\Gamma u, \nabla^\Gamma u \rangle \, d\mu < \infty$$

and such that

(1.3) $\quad \dfrac{u \circ \psi_t^* - u}{t} \longrightarrow \langle \nabla^\Gamma u, V \rangle \qquad \text{as } t \to 0 \text{ in } L^2(\mu \circ (\psi_s^*)^{-1}), \text{ for all } s \in \mathbb{R},$

whenever $V \in V_0(X)$ with flow $(\psi_t)_{t \in \mathbb{R}}$. By $\mathbb{F}^{(c)}$ we denote the subset of all continuous elements of $\mathbb{F}$. Clearly $\mathcal{F}C_b^\infty \subset \mathbb{F}^{(c)}$ by Assumption 1.1 (b).

**Proposition 1.4:** *Assume that $\mu$ satisfies Assumption 1.1.*
(i) *$(\mathcal{E}_\mu^\Gamma, \mathcal{F}C_b^\infty)$, $(\mathcal{E}_\mu^\Gamma, \mathbb{F})$, and $(\mathcal{E}_\mu^\Gamma, \mathbb{F}^{(c)})$ are closable and their closures, denoted by $(\mathcal{E}_\mu^\Gamma, \mathcal{F}_0)$, $(\mathcal{E}_\mu^\Gamma, \mathcal{F})$, and $(\mathcal{E}_\mu^\Gamma, \mathcal{F}^{(c)})$ respectively, are Dirichlet forms. Clearly, $\mathcal{F}_0 \subset \mathcal{F}^{(c)} \subset \mathcal{F}$.*
(ii) *If $u \in L^2(\mu)$ is $\rho$-Lipschitz continuous, then $u \in \mathcal{F}$.*
(iii) *Suppose that $\mu = \int \pi_{s \cdot \sigma} \lambda(ds)$ is as in Example 1.2 and that the $X$-valued Brownian motion having the drift $\nabla^X \log(d\sigma/dm)$ is conservative. Then $\mathcal{F}_0 = \mathcal{F}^{(c)} = \mathcal{F}$.*

(iv) *For each $u \in \mathcal{F}$, there is a section $\nabla^\Gamma u$ of the tangent bundle $T\Gamma_X$ such that $\mathcal{E}_\mu^\Gamma(u,u) = \int \langle \nabla^\Gamma u(\gamma), \nabla^\Gamma u(\gamma) \rangle_\gamma \, \mu(d\gamma)$ I.e., our Dirichlet form admits a* carré du champs operator *given by $\mathbf{\Gamma}_\mathcal{E}(u,v)(\gamma) = \langle \nabla^\Gamma u(\gamma), \nabla^\Gamma v(\gamma) \rangle_\gamma$.*

Proposition 1.4 will be proved in Section 6 below. A priori it is not clear whether one of the identities in (iii) transfer to more general cases than mixed Poisson measures. An investigation in this direction will be the subject of a future work. Our next result first gives a partial converse to Theorem 1.3. Then we show in (ii) that $\rho$ is in fact the largest metric which yields the assertion of this theorem.

**Theorem 1.5:** *Suppose that $\mu$ satisfies Assumption 1.1.*
(i) *If $u \in \mathcal{F}$ satisfies $\mathbf{\Gamma}_\mathcal{E}(u,u) \leq C^2$ $\mu$-a.e. and if $u$ has a $\rho$-continuous $\mu$-version, then there exists a $\mu$-measurable $\mu$-version $\widetilde{u}$ which is $\rho$-Lipschitz continuous and satisfies $\mathrm{Lip}(\widetilde{u}) \leq C$.*
(ii) *$\rho$ is the intrinsic metric of the Dirichlet forms $(\mathcal{E}_\mu^\Gamma, \mathcal{F})$ and $(\mathcal{E}_\mu^\Gamma, \mathcal{F}^{(c)})$, i.e.,*
$\rho(\gamma, \omega) =$

$$
\begin{aligned}
(1.4) \quad &= \sup \left\{ u(\gamma) - u(\omega) \,\Big|\, u \in \mathcal{F} \cap C(\Gamma_X) \text{ and } \mathbf{\Gamma}_\mathcal{E}(u,u) \leq 1 \ \mu\text{-a.e. on } \Gamma_X \right\} \\
&= \sup \left\{ u(\gamma) - u(\omega) \,\Big|\, u \in \mathbb{F}^{(c)} \text{ and } \mathbf{\Gamma}_\mathcal{E}(u,u) \leq 1 \ \mu\text{-a.e. on } \Gamma_X \right\}.
\end{aligned}
$$

**Remark:** We believe that in (i) the assumption that $u$ has a $\rho$-continuous version can be dropped, at least if $\mu$ is a mixed Poisson measure. Theorem 1.5 (ii) states in particular that continuous functions $u \in \mathbb{F}^{(c)}$ with $\mathbf{\Gamma}_\mathcal{E}(u,u)$ bounded are already $\rho$-Lipschitz continuous. However, if $X$ is not compact an arbitrary $\rho$-Lipschitz continuous $u$ will in general have uncountably many $\rho$-Lipschitz continuous $\mu$-versions with arbitrarily large Lipschitz constant, which can be seen by modifying $u$ on a single fiber $\{\rho(\gamma, \cdot) < \infty\}$ having measure 0. Therefore, it would not make sense to replace $\mathbb{F}^{(c)}$ or $\mathcal{F} \cap C(\Gamma_X)$ in (1.4) by a larger class of not necessarily continuous functions. On the other hand, it could be useful to know whether $\mathbb{F}^{(c)}$ can be replaced by the smaller set $\mathcal{F}C_b^\infty$. Note that (1.4) is reminiscent of the well-known Kantorovich-Rubinstein theorem for the $L^1$-Wasserstein metric between probability measures. One might guess that a similar variational formula as (1.4) holds for the classical $L^2$-Wasserstein distance on the space of probabilities. It seems that in this case only variational characterizations involving non-symmetric expressions are known to date (cf. Dudley (1989)). If $X$ is compact, then $\rho$ metrizes the vague topology on $\Gamma_X$. Hence the well-known results on intrinsic distances of regular Dirichlet forms apply (see e.g. Sturm (1995)). However, if $X$ is not compact the situation changes completely, and one is reminded of the Cameron-Martin distance on path space.

## 2. Applications to Ruelle-type Gibbs measures

Let us recall the terminology used for Ruelle measures. Suppose $X = \mathbb{R}^d$, and $\phi : \mathbb{R}^d \to \mathbb{R} \cup \{+\infty\}$ is a *pair potential*, i.e., $\phi$ is measurable and satisfies $\phi(x) = \phi(-x)$. If $\Lambda \subset \mathbb{R}^d$ is open bounded and non-empty the conditional energy $E_\Lambda^\phi : \Gamma_{\mathbb{R}^d} \to \mathbb{R} \cup \{+\infty\}$ is defined as

$$E_\Lambda^\phi(\gamma) := \int_\Lambda \int_{\mathbb{R}^d} \phi(x-y) \mathrm{I}_{\{x \neq y\}} \gamma(dx)\gamma(dy)$$

if

$$\int_\Lambda \int_{\mathbb{R}^d} |\phi(x-y)| \mathrm{I}_{\{x \neq y\}} \gamma(dx)\gamma(dy) < \infty,$$

and $E_\Lambda^\phi(\gamma) := +\infty$ otherwise. For $r = (r^1, \ldots, r^d) \in \mathbb{Z}^d$, let $Q_r$ denote the cube

$$Q_r := \left\{ x \in \mathbb{R}^d \,\Big|\, r^i - \frac{1}{2} \leq x^i < r^i + \frac{1}{2} \right\},$$

and, for $N \in \mathbb{N}$ define $\Lambda_N := [-N, N]^d$. Then $\phi$ is called *superstable*, if there are constants $A > 0$ and $B \geq 0$ such that

$$E_{\Lambda_N}^\phi(\gamma_{\Lambda_N}) \geq \sum_{r \in \mathbb{Z}^d} \left[ A \gamma_{\Lambda_N}(Q_r)^2 - B \gamma_{\Lambda_N}(Q_r) \right], \qquad \text{for all } \gamma \in \Gamma_{\mathbb{R}^d} \text{ and } N \in \mathbb{N}.$$

Recall that $\gamma_\Lambda$ denotes the restriction of $\gamma$ to $\Lambda$. Let $|\cdot|_\infty$ denote the maximum norm on $\mathbb{R}^d$. $\phi$ is called *lower regular* if there exists a decreasing positive function $a : \mathbb{N} \to [0, \infty)$ such that $\sum_{r \in \mathbb{Z}^d} a(|r|_\infty) < \infty$ and, for any $\Lambda'$, $\Lambda''$ which are finite unions of cubes $Q_r$ and disjoint,

$$W(\gamma_{\Lambda'}|\gamma_{\Lambda''}) := \int_{\Lambda'} \int_{\Lambda''} \phi(x-y)\,\gamma(dx)\,\gamma(dy)$$
$$\geq -\sum_{r',r'' \in \mathbb{Z}^d} a(|r'-r''|_\infty)\gamma_{\Lambda'}(Q_{r'})\gamma_{\Lambda''}(Q_{r''}),$$

for all $\gamma \in \Gamma_{\mathbb{R}^d}$. See Ruelle (1970) and the references therein for a discussion of these conditions. Now let $m$ denote Lebesgue measure on $\mathbb{R}^d$ and let $z > 0$ be fixed. Then let $Z_\Lambda^\phi : \Gamma_{\mathbb{R}^d} \to [0, \infty]$ be the partition function

$$Z_\Lambda^\phi(\gamma) = \int \exp\left(-E_\Lambda^\phi(\gamma_{\Lambda^C} + \omega_\Lambda)\right) \pi_{z \cdot m}(d\omega),$$

where $\Lambda^C := \mathbb{R}^d \setminus \Lambda$. For $B \subset \Gamma_{\mathbb{R}^d}$ measurable, we define

$$\Pi_\Lambda^\phi(\gamma, B) = \mathrm{I}_{\{Z_\Lambda^\phi < \infty\}}(\gamma) \frac{1}{Z_\Lambda^\phi(\gamma)} \int \mathrm{I}_B(\gamma_{\Lambda^C} + \omega_\Lambda) \exp\left(-E_\Lambda^\phi(\gamma_{\Lambda^C} + \omega_\Lambda)\right) \pi_{z \cdot m}(d\omega).$$

The system $\Pi_\Lambda^\phi$, $\Lambda \subset \mathbb{R}^d$ open and bounded, is a specification and $\mu$ is a Gibbs measure with respect to $\Pi_\Lambda^\phi$ if it satisfies the equilibrium equations

$$\mu \Pi_\Lambda^\phi = \mu, \qquad \text{for all } \Lambda.$$

Because of the way our specification was defined, $\mu$ is also called *grandcanonical Gibbs measure associated with* $\phi$. Such a measure $\mu$ is called *tempered* if it is supported by $S_\infty := \bigcup_{n=1}^\infty S_n$, where

$$S_n := \Big\{ \gamma \in \Gamma_{\mathbb{R}^d} \ \Big| \ \forall\, N \in \mathbb{N}, \sum_{r \in \Lambda_N \cap \mathbb{Z}^d} \gamma(Q_r)^2 \leq n^2 (2N+1)^d \Big\}.$$

According to Section 5 of Ruelle (1970), the set of tempered grandcanonical Gibbs measures is non-empty, provided the potential $\phi$ is superstable, lower regular, and satisfies the following integrability condition

$$(2.1) \qquad \int_{\mathbb{R}^d} \big|1 - e^{-\phi(x)}\big|\, dx < +\infty.$$

The following differentiability condition on $\phi$ was introduced in Albeverio, Kondratiev and Röckner (1997b).

(2.2) $e^{-\phi}$ is weakly differentiable on $\mathbb{R}^d$, $\phi$ is weakly differentiable on $\mathbb{R}^d \backslash \{0\}$, and the weak gradient $\nabla \phi$ (which is a locally $m$-integrable function on $\mathbb{R}^d \backslash \{0\}$) considered as an $m$-a.e. defined function on $\mathbb{R}^d$ satisfies $\nabla \phi \in L^1(\mathbb{R}^d, e^{-\phi} dm) \cap L^2(\mathbb{R}^d, e^{-\phi} dm)$.

**Proposition 2.1:** *Suppose $\phi : \mathbb{R}^d \to \mathbb{R}^d \cup \{+\infty\}$ is a superstable lower regular pair potential with compact support which in addition satisfies (2.1) and (2.2) and which is bounded on any set $\{x\,|\,|x| > r\}$, for all $r > 0$. Then every tempered grandcanonical Gibbs measure associated with $\phi$ (in short: Ruelle measure) satisfies Assumption 1.1.*

Together with Proposition 5.6 below this result immediately yields the following corollary.

**Corollary 2.2:** *If $\phi$ is as in Proposition 2.1, then every Ruelle measure associated with $\phi$ has full topological support on $\Gamma_X$.*

**Proof of Proposition 2.1:** Suppose first that $\phi$ is superstable, lower regular, and satisfies (2.1). Then Corollary 5.3 of Ruelle (1970) states that, for any tempered grandcanonical Gibbs measure $\mu$ and any bounded open set $\Lambda \subset \mathbb{R}^d$, there exists $\sigma_\Lambda : \Gamma_{\mathbb{R}^d} \to [0, \infty)$ such that, for any measurable function $F \geq 0$,

$$\int F(\gamma_\Lambda) \, \mu(d\gamma) = \int \sigma_\Lambda(\gamma) F(\gamma_\Lambda) \, \pi_m(d\gamma).$$

Moreover, there are constants $c > 0$ and $d \in \mathbb{R}$ such that

$$\sigma_\Lambda(\gamma) \leq \exp\Big(m(\Lambda) + \sum_{r \in \mathbb{Z}^d} \big[ - c\gamma(\Lambda \cap Q_r)^2 + d\gamma(\Lambda \cap Q_r) \big]\Big),$$

for $\pi_m$-a.e. $\gamma \in \Gamma_{\mathbb{R}^d}$. In particular, the density $\sigma_\Lambda$ is bounded above by a constant. It follows immediately from this result that (a) and (b) of Assumption 1.1 are satisfied for such $\mu$. Condition (c) follows from Lemma 2.3 below. Finally, Albeverio, Kondratiev and Röckner (1997b) show that the quasi-invariance and the integration by parts formula of Assumption 1.1 (d) and (e) respectively hold under our assumptions (cf. Lemma 4.2, Theorem 4.3, and Section 5.1 of Albeverio, Kondratiev and Röckner (1997b)). Moreover, if $V \in V_0(\mathbb{R}^d)$ has support contained in $\Lambda$ and generates the flow $(\psi_t)_{t \in \mathbb{R}}$, then

$$\frac{d\mu \circ (\psi_t^*)^{-1}}{d\mu}(\gamma) = \exp\big[E_\Lambda^\phi(\gamma) - E_\Lambda^\phi(\psi_{-t}^*\gamma)\big] \frac{d\pi_{z \cdot m} \circ (\psi_t^*)^{-1}}{d\pi_{z \cdot m}}(\gamma),$$

and the positivity condition in (d) holds under our assumptions on $\phi$. Thus Proposition 2.1 is proved. $\square$

We owe the following lemma to B. Schmuland.

**Lemma 2.3:** *Let $\mu$ be grandcanonical Gibbs measure with respect to $\phi$. Then*

$$\mu\big(\{\gamma \,|\, \gamma(\mathbb{R}^d) < \infty\}\big) = 0.$$

**Proof:** Let 0 denote the empty configuration. First suppose that $Z_\Lambda^\phi(0) < \infty$, for all bounded open $\Lambda \subset \mathbb{R}^d$. Then

$$\mu(\{0\}) = \int \frac{1}{Z_\Lambda^\phi(\gamma)} \int_{\{\omega_\Lambda + \gamma_{\Lambda^c} = 0\}} \exp\big[-E_\Lambda^\phi(\omega_\Lambda + \gamma_{\Lambda^c})\big] \, \pi_{z \cdot m}(d\omega) \, \mu(d\gamma)$$

(2.3)
$$= \int_{\{\gamma_{\Lambda^c} = 0\}} \frac{1}{Z_\Lambda^\phi(\gamma)} \exp[-z \cdot m(\Lambda)] \, \mu(d\gamma)$$

$$= \frac{1}{Z_\Lambda^\phi(0)} \exp[-z \cdot m(\Lambda)] \cdot \mu(\gamma_{\Lambda^c} = 0).$$

Now

$$Z_\Lambda^\phi(0) = \int \exp\big[-E_\Lambda^\phi(\omega_\Lambda)\big]\,\pi_{z\cdot m}(d\omega) \geq \exp[-z\cdot m(\Lambda)]\cdot(1+z\cdot m(\Lambda)),$$

because $E_\Lambda^\phi(\omega_\Lambda) = 0$ if $\omega(\Lambda)$ equals 0 or 1. Together with (2.3) this yields

$$\mu(\{0\}) \leq \frac{1}{1+z\cdot m(\Lambda)} \longrightarrow 0 \quad \text{as } \Lambda \uparrow \mathbb{R}^d.$$

Plugging this back into (2.3) gives $\mu(\gamma_{\Lambda^C} = 0) = 0$, for all $\Lambda$, which in turn implies that $\mu(\gamma(\mathbb{R}^d) < \infty) = 0$.

Next consider the case where $Z_\Lambda^\phi(0) = \infty$, for large $\Lambda$. Then

$$\mu\big(\{\gamma \,|\, \gamma(\mathbb{R}^d) < \infty\}\big) = \int \Pi_\Lambda^\phi\big(\gamma, \{\omega \,|\, \gamma(\Lambda^C) + \omega(\Lambda) < \infty\}\big)\,\mu(d\gamma)$$

$$= \int_{\{\gamma(\mathbb{R}^d)<\infty\}} \Pi_\Lambda^\phi\big(\gamma, \{\omega \,|\, \gamma(\Lambda^C) + \omega(\Lambda) < \infty\}\big)\,\mu(d\gamma).$$

But if $\gamma(\mathbb{R}^d) < \infty$, then $\Pi_\Lambda^\phi(\gamma, \cdot) = \Pi_\Lambda^\phi(0, \cdot)$, for $\Lambda$ large. Also $\Pi_\Lambda^\phi(0, \cdot)$ is the zero measure, for $\Lambda$ so large that $Z_\Lambda^\phi(0) = \infty$. This proves the lemma. □

## 3. Application: Potential theory on configuration space

Let $\mu$ satisfy Assumption 1.1, and define, for $A \subset \Gamma_X$,

$$\rho_A(\gamma) := \inf\big\{\rho(\omega,\gamma)\,\big|\,\omega \in A\big\}.$$

It will be shown in Lemma 4.1 below that $\rho_A$ is a measurable function if $A$ is closed.

**Proposition 3.1:** *If $K \subset \Gamma_X$ is compact and $c \geq 0$, then $c \wedge \rho_K$ is an $\mathcal{E}_\mu^\Gamma$-quasi continuous function in $\mathcal{F}^{(c)}$.*

**Proof:** For $\omega \in \Gamma_X$ and $r > 0$, let $\rho_{\omega,r}$ denote the function defined in Lemma 4.2 below, and let $\rho_{K,r}(\gamma) = \inf\{\rho_{\omega,r}(\gamma)\,|\,\omega \in K\}$. Let $\omega_\Lambda$ denote the restriction of a configuration $\omega$ to $\Lambda \subset X$. If $F_r = \big\{\widetilde{\omega}\,\big|\,\exists\omega \in K \text{ such that } \widetilde{\omega}_{B_r} = \omega_{B_r}\big\}$, then $F_r$ is closed, and $\rho_{K,r} = \rho_{F_r}$. Hence $\rho_{K,r}$ is lower semi-continuous by Lemma 4.1 (vii). However, $\rho_{K,r}$ is also upper semi-continuous as infimum over the continuous functions $\rho_{\omega,r}$. Hence $\rho_{K,r}$ is continuous.

Let us now show that $\lim_{r\uparrow\infty} \rho_{K,r}(\gamma) = \rho_K(\gamma)$, for all $\gamma \in \Gamma_X$. We first note that the limit exists since $r \mapsto \rho_{K,r}$ is increasing. Furthermore, by Lemma 4.1 there is $\omega_\gamma \in K$ such that

$$\rho_K(\gamma) = \rho(\omega_\gamma,\gamma) = \lim_{r\uparrow\infty} \rho_{\omega_\gamma,r}(\gamma) \geq \lim_{r\uparrow\infty} \rho_{K,r}(\gamma),$$

where we have used Lemma 4.2 (iii) for the second identity. To prove that also $\rho_K(\gamma) \leq \lim_r \rho_{K,r}(\gamma)$ choose $\alpha < \rho_K(\gamma)$. By Lemma 4.2 (iv), the set $U_r := \{\omega \,|\, \rho_{\omega,r}(\gamma) > \alpha\}$ is open. In addition, Lemma 4.2 (iii) implies that, for any $\omega \in K$, there is some $r_\omega \in (0, \infty)$ such that $\rho_{\omega, r_\omega}(\gamma) > \alpha$. Thus $\{U_r \,|\, r > 0\}$ covers $K$, and hence there must be some $r^{(\alpha)} < \infty$ with $K \subset U_{r^{(\alpha)}}$. In other words, $\rho_{K, r^{(\alpha)}}(\gamma) \geq \alpha$. Hence we conclude that $\rho_K(\gamma) \leq \lim_r \rho_{K,r}(\gamma)$.

Next, Theorem 1.3 and Lemma 4.1 (viii) yield that $u_r := c \wedge \rho_{K,r} \in \mathbb{F}^{(c)} \subset \mathcal{F}^{(c)}$ and that $\mathcal{E}^\Gamma_\mu(u_r, u_r) \leq 1$. But this implies the assertion (cf. Ma and Röckner (1992), Section III.3). $\square$

**Corollary 3.2:** If $\mu(A) = 1$, then $\{\rho_A > 0\}$ is exceptional for $(\mathcal{E}^\Gamma_\mu, \mathcal{F}^{(c)})$.

**Proof:** By inner regularity there are compact sets $K_1 \subset K_2 \subset \cdots \subset A$ such that $\mu(K_n) \uparrow 1$. Consider the functions $u_n(\gamma) := \rho_{K_n}(\gamma) \wedge 1$. Then $u_n \in \mathcal{F}^{(c)}$, and $\mathcal{E}^\Gamma_\mu(u_n, u_n) \leq 1$ by Proposition 3.1, Theorem 1.3, and Lemma 4.1 (viii) below. Therefore $u := \lim_n u_n$ is non-negative and $\mathcal{E}^\Gamma_\mu$-quasi continuous by Proposition 3.1 and standard arguments (see e.g. Ma and Röckner (1992) Section III.3). In addition $u = 0$ on $\bigcup_n K_n$, and hence $u = 0$ $\mu$-a.s. Hence $u = 0$ even quasi everywhere by Proposition III.3.9 of Ma and Röckner (1992). But $\{u = 0\} \subset \{\rho_A > 0\}$. $\square$

Let us now look at some applications of Corollary 3.2:

**Example 3.3:**
(i) If $A = \{\gamma \,|\, \gamma(X) = \infty\}$ one sees immediately that $A = \{\rho_A < \infty\}$. Hence $\mu(A) = 1$ implies that $A^C$ is exceptional, whereas $\mu(A) = 0$ implies that $A$ is exceptional. this result has first been proved by Byron Schmuland (private communication).

(ii) For the next application, suppose $X = \mathbb{R}^d$. Consider the shift transformation $\theta_x$, which is defined by $\theta_x \gamma = \delta_x * \gamma$. Note that $\rho(\omega, \gamma) < \infty$ implies that

$$|u(\theta_x \omega) - u(\theta_x \gamma)| \longrightarrow 0 \quad \text{as } |x| \to \infty,$$

for all functions $u$ which are uniformly continuous with respect to a metric $\delta$ for the vague topology on $\Gamma_{\mathbb{R}^d}$ having the form $\delta(\omega, \gamma) = \sum_{i=1}^{\infty} 2^{-i} \big| \int f_i \, d\omega - \int f_i \, d\gamma \big| \wedge 1$ with convergence determining $f_i \in C_0(\mathbb{R}^d)$. Now let

$$A_\mu = \left\{ \gamma \in \Gamma_{\mathbb{R}^d} \ \Big| \ \frac{1}{|V_n|} \int_{V_n} \delta_{\theta_x \gamma} \, dx \longrightarrow \mu \text{ as } n \uparrow \infty \right\},$$

where $V_n$ is the box $[-n, n]^d$ and convergence is supposed to hold in the weak sense. Then again $A_\mu = \{\rho_{A_\mu} < \infty\}$. But if $\mu$ is ergodic with respect to $\theta_x$ the spatial ergodic theorem of Nguyen and Zessin (1979) implies $\mu(A_\mu) = 1$. Thus $A_\mu^C$ is exceptional by our corollary. If $\mu$ is not ergodic one can use its ergodic decomposition to get an analogous result.

(iii) A similar reasoning as above applies to the strong law of large numbers.

The above examples exhibit an interesting relation between the tail structure of $\Gamma_{\mathbb{R}^d}$ and the $\sigma$-field of all events $A$ with the property that

$$A = \{\rho_A < \infty\}.$$

This relation and its application to Gibbs measures will be exploited in a future work.

Next we present a short proof for the quasi-regularity of our Dirichlet form $(\mathcal{E}_\mu^\Gamma, \mathcal{F}^{(c)})$. This property implies in particular that the form is associated with a diffusion process having $\mu$ as symmetrizing measure. Hence our exceptional sets above can be interpreted as polar sets for this diffusion process. Note that an slightly stronger result, the quasi-regularity of $(\mathcal{E}_\mu^\Gamma, \mathcal{F}_0)$, has been proved in Ma and Röckner (1997). We refer to Ma and Röckner (1992), Chapter IV, for the terminology below.

**Corollary 3.4:** $(\mathcal{E}_\mu^\Gamma, \mathcal{F}^{(c)})$ is quasi-regular.

**Proof:** It suffices to show that Cap is tight. To this end, let $K_1 \subset K_2 \subset \cdots \subset \Gamma_X$ be compact with $\mu(K_n) \to 1$. Then note that the sets $\{\rho_{K_n} \leq 1/2\}$ are also compact by Lemma 4.1 (vii) below. But by standard arguments there exists a subsequence $(n_m)_{m \in \mathbb{N}}$ such that $u_N := N^{-1} \sum_{m=1}^N \rho_{K_{n_m}} \wedge 1 \to 0$ wrt $\mathcal{E}_\mu^\Gamma + (\cdot, \cdot)_{L^2(\mu)}$. Moreover, $u_N \geq 1/2$ on $\{\rho_{K_{n_N}} \geq 1/2\}$. Hence

$$\mathrm{Cap}(\rho_{K_{n_N}} > 1/2) \leq \mathrm{Cap}(u_N > 1/2) \leq \mathcal{E}_\mu^\Gamma(u_N, u_N) + \int u_N^2 \, d\mu \longrightarrow 0$$

by Proposition III.3.4 of Ma and Röckner (1992).  □

## 4. Topological properties of $\rho$

In this section, we collect some preliminary results of topological kind concerning our metric $\rho$. Recall that $\Gamma_X$ is always endowed with the vague topology.

**Lemma 4.1:** Let $\pi_i : X \times X \to X$, $i = 1, 2$ denote the projection on the i-th coordinate
  (i) The mapping $\Gamma_{X \times X} \ni \eta \mapsto \int d(x, y)^2 \, \eta(dx, dy)$ is lower semi-continuous.
 (ii) If $K \subset \Gamma_X$ is compact and $i \in \{1, 2\}$, then the set $\{\eta \mid \pi_i^* \eta \in K\}$ is relatively compact in $\Gamma_{X \times X}$.

(iii) For every $\alpha > 0$, $i = 1, 2$, the projection map $\pi_i^* : \Gamma_{X \times X} \to \Gamma_X$ restricted to the closed set $G_\alpha := \{\eta \in \Gamma_{X \times X} \mid \int d(x, y)^2 \, \eta(dx, dy) \leq \alpha^2\}$ is continuous.

(iv) Suppose $\gamma, \omega \in \Gamma_X$ have finite $\rho$-distance. Then there is at least one $\eta^* \in \Gamma_{\gamma, \omega}$ such that $\rho(\gamma, \omega) = \left(\int d(x, y)^2 \, \eta^*(dx, dy)\right)^{1/2}$

(v) The map $G_\alpha \ni \eta \mapsto (\pi_1^* \eta, \pi_2^* \eta) \in \Gamma_X \times \Gamma_X$ is closed (i.e., $\{(\pi_1^* \eta, \pi_2^* \eta) \mid \eta \in F \cap G_\alpha\}$ is closed in $\Gamma_X \times \Gamma_X$, for all closed $F \subset \Gamma_{X \times X}$). In particular, $\pi_1^*, \pi_2^* : G_\alpha \to \Gamma_X$ are both closed.

(vi) $\rho$ is lower semi-continuous on $\Gamma_X \times \Gamma_X$.

(vii) Let $A \subset \Gamma_X$, $A$ closed. Then $\rho_A$ is lower semi-continuous (hence measurable) on $\Gamma_X$. If $A$ is compact, then $\{\rho_A \leq \alpha\}$ is compact, for all $\alpha \geq 0$. In particular, closed $\rho$-balls are compact.

(viii) $\rho_A \wedge c$ is $\rho$-Lipschitz continuous with Lipschitz constant $\leq 1$, for all $A \subset \Gamma_X$ and $c \geq 0$.

**Proof:** (i) is trivial. For the proof of (ii), let $C$ denote the set under consideration. Then $C$ is relatively compact if and only if $\sup_{\eta \in C} \eta(F) < \infty$, for all compact sets $F \subset X \times X$. But, for such $F$,

$$\sup_{\eta \in C} \eta(F) \leq \sup_{\eta \in C} \pi_i^* \eta(\pi_i(F)) \leq \sup_{\gamma \in K} \gamma(\pi_i(F)) < \infty$$

by compactness of $K$.

To prove (iii), let $f \in C_0(X)$, and choose $g \in C_0(X)$ so that $g \equiv 1$ on $(\operatorname{supp} f)^\alpha := \{d_{\operatorname{supp} f} \leq \alpha\}$. Then no $\eta \in G_\alpha$ charges a point $(x, y)$ where $x \in \operatorname{supp} f$ and $y \notin (\operatorname{supp} f)^\alpha$. Hence, for all $\eta \in G_\alpha$,

$$\int f(x) \, \pi_1^* \eta(dx) = \int f(x) \, \eta(dx, dy) = \int f(x) g(y) \, \eta(dx, dy).$$

In particular, $\pi_1^* \eta \in \Gamma_X$ and it follows that $\pi_1^*$, restricted to $G_\alpha$, is continuous.

(iv) follows from (i), (ii), and (iii). For the proof of (v), let $F \subset \Gamma_{X \times X}$, $F$ closed, and $\eta_k \in F \cap G_\alpha$, $k \in \mathbb{N}$, such that $\pi_i^* \eta_k \to \gamma_i \in \Gamma_X$, $i = 1, 2$. Then

$$\eta_k \in \left\{\eta \in F \cap G_\alpha \,\Big|\, \pi_1^* \eta \in \{\pi_1^* \eta_k \mid k \in \mathbb{N}\} \cup \{\gamma_1\}\right\} =: C,$$

which is a compact set by (i) and (ii). Hence selecting a subsequence if necessary we may assume that $\eta_k \to \eta' \in C \subset F \cap G_\alpha$. Consequently, by continuity $\gamma_i = \lim_k \pi_i^* \eta_k = \pi_i^* \eta$. Hence $(\gamma_1, \gamma_2) \in \{(\pi_1^* \eta, \pi_2^* \eta) \mid \eta \in F \cap G_\alpha\}$.

(vi). Let $\alpha \geq 0$. Then by (iv)

$$\left\{(\gamma, \omega) \in \Gamma_X \times \Gamma_X \,\Big|\, \rho(\omega, \gamma) \leq \alpha\right\} = \left\{(\pi_1^* \eta, \pi_2^* \eta) \,\Big|\, \eta \in G_\alpha\right\}.$$

But by (v) the latter is a closed set.

(vii). Let $\alpha \geq 0$ and $\gamma \in \Gamma_X$. Then

$$\overline{B}_\alpha^\rho(\gamma) := \{\omega \in \Gamma_X \mid \rho(\gamma,\omega) \leq \alpha\} = \{\pi_2^*\eta \mid \eta \in G_\alpha,\, \pi_1^*\eta = \gamma\}$$

is compact as the continuous image of $\{\eta \in G_\alpha \mid \pi_1^*\eta = \gamma\}$. In particular, $\omega \mapsto \rho(\gamma,\omega)$ is lower semi-continuous, and if $r := \rho_A(\gamma) < \infty$ there exists $\omega_\gamma \in C := \overline{B}_{r+1}^\rho(\gamma) \cap A$ such that $\rho(\gamma,\omega_\gamma) = \rho_C(\gamma) = \rho_A(\gamma)$. Hence

$$\{\rho_A \leq \alpha\} = \{\pi_1^*\eta \mid \eta \in G_\alpha,\, \pi_2^*\eta \in A\}.$$

But the latter set is closed by (v), and even compact if $A$ is by (ii) and (iii).

(viii) is trivial, and the lemma is proved. $\square$

The next Lemma will in particular imply the quasi-continuity of $\rho(\omega,\cdot)$, for fixed $\omega$. Note that its Assertion (i) needs that $X$ is connected.

**Lemma 4.2:** *Let $B_r$ denote the open geodesic ball of radius $r > 0$ centered in some fixed point of $X$, and recall that $\gamma_{B_r}$ denotes the restriction of $\gamma \in \Gamma_X$ to $B_r$. For $\omega \in \Gamma_X$, we define the closed set $A_{\omega,r} := \{\gamma \in \Gamma_X \mid \gamma_{B_r} = \omega_{B_r}\}$, and let $\rho_{\omega,r} := \rho_{A_{\omega,r}}$.*

*(i) If $B_r \neq X$, then $\rho_{\omega,r}(\gamma) < \infty$ if and only if $\gamma(X) \geq \omega(B_r)$.*

*(ii) $\rho_{\omega,r}$ is a continuous function from $\Gamma_X$ to $[0,\infty]$.*

*(iii) $\rho_{\omega,r}(\gamma) \nearrow \rho(\omega,\gamma)$ as $r \uparrow \infty$, for all pairs $\omega, \gamma \in \Gamma_X$.*

*(iv) $\omega \mapsto \rho_{\omega,r}(\gamma)$ is lower semi-continuous if $\gamma$ is fixed.*

**Proof:** (i). First fix $\omega \in \Gamma_X$ and write $\omega_{B_r}$ as $\sum_{i=1}^n \delta_{y^i}$. If $\gamma \in \Gamma_X$ is given write it as $\sum_{i \in I} \delta_{z^i}$, for some index set $I \subset \mathbb{N}$. If $|I| < n$ it is clear that $\rho_{\omega,r}(\gamma) = \infty$. Therefore we can assume $|I| \geq n$ in the sequel. Then we can find $n$ points $x^1, \ldots, x^n$ such that $\gamma' := \gamma - \sum_{i=1}^n \delta_{x^i}$ is a non-negative measure. Then we can write $\gamma'_{B_r}$ as $\sum_{i=1}^m \delta_{x^{i+n}}$, for some $m \geq 0$ and $x^{n+1}, \ldots, x^{n+m} \in B_r$. For $i = 1, \ldots, m$, we then pick an $y^{n+i} \in \partial B_r$ which realizes the Riemannian distance of $x^{n+i}$ to the boundary $\partial B_r$ of $B_r$. If we now define an element $\widetilde{\omega}$ of $\Gamma_X$ by $\widetilde{\omega} := \gamma'_{B_r^C} + \sum_{i=1}^{n+m} \delta_{y^i}$ it is clear that $\widetilde{\omega}_{B_r} = \omega_{B_r}$ and that $\rho(\widetilde{\omega},\gamma)^2 \leq r^2 m + \sum_{i=1}^n d(x^i, y^i)^2 < \infty$.

(ii). In view of Lemma 4.1 (vii) it suffices to show upper semicontinuity of $\gamma \mapsto \rho_{\omega,r}(\gamma) = \rho_{A_{\omega,r}}(\gamma)$. To this end, suppose we are given a sequence $(\gamma_k)_{k \in \mathbb{N}} \subset \Gamma_X$ with $\gamma_k \to \gamma$. If $\rho_{\omega,r}(\gamma) = \infty$ we are done. Thus assume $\rho_{\omega,r}(\gamma) < \infty$ in the sequel. Note that, for any $\omega' \in A_{\omega,r}$ and $\eta \in \Gamma_{X \times X}$ so that

$$\pi_1^*\eta = \gamma, \quad (\pi_2^*\eta)_{B_r} = \omega_{B_r}, \quad \rho(\gamma,\omega')^2 = \int d(x,y)^2 \, \eta(dx,dy),$$

we can construct a new $\eta'$ by replacing any $(x,y) \in \operatorname{supp} \eta$ with $x, y \in B_r^C$ by $(x,x)$. Then

$$\pi_1^* \eta' = \gamma, \quad (\pi_2^* \eta')_{B_r} = \omega_{B_r}, \quad \text{but } \int d(x,y)^2 \eta'(dx, dy) \leq \int d(x,y)^2 \eta(dx, dy).$$

Since $A_{\omega,r}$ is closed, we can hence find some $\omega^* \in A_{\omega,r}$ and $\eta^* \in \Gamma_{\gamma,\omega^*}$ such that $\rho_{\omega,r}(\gamma) = \rho(\omega^*, \gamma)$, and such that $\eta^*$ is optimal in the sense of 4.1 (iv) and has the form $\eta^* = \sum_{i=1}^{N} \delta_{(x^i, y^i)} + \sum \delta_{(x,x)}$, for $(x^i, y^i) \in B_r \times \Gamma_X \cup \Gamma_X \times B_r$. Let $\gamma' = \gamma - \sum_{i=1}^{N} \delta_{x^i}$, and write $\gamma'_{\partial B_r}$ as $\sum_{i=1}^{m} \delta_{x^{N+i}}$. Then there exist $x_k^1, \ldots x_k^{N+m}$ in $\gamma_k$ such that $x_k^i \to x^i$. Define, for $k \in \mathbb{N}$,

$$\eta_k := \sum_{i=1}^{N} \delta_{(x_k^i, y^i)} + \sum_{i=N+1}^{N+m} \delta_{(x_k^i, x^i)} + \sum \delta_{(x_k, x_k)},$$

where third sum is determined by $\pi_1^* \eta_k = \gamma_k$. Then, for large $k$, $(\pi_2^* \eta_k)_{B_r} = \omega_{B_r}$, since $(\gamma_k - \sum_{i=1}^{N+m} \delta_{x_k^i})(\overline{B}_r) = 0$ eventually. Therefore,

$$\limsup_{k \uparrow \infty} \rho_{\omega,r}(\gamma_k)^2 \leq \limsup_{k \uparrow \infty} \Big( \sum_{i=1}^{N} d(x_k^i, y^i)^2 + \sum_{i=N+1}^{N+m} d(x_k^i, x^i)^2 \Big)$$

$$= \sum_{i=1}^{N} d(x^i, y^i)^2 = \rho_{\omega,r}(\gamma)^2.$$

(iii). It suffices to show that $\alpha := \sup_r \rho_{\omega,r}(\gamma) \geq \rho(\omega, \gamma)$. If $\alpha = \infty$ we are done. Otherwise we know from the proof of (ii) that there are $\omega_r^*$ such that $(\omega_r^*)_{B_r} = \omega_{B_r}$ and $\rho_{\omega,r}(\gamma) = \rho(\omega_r^*, \gamma)$. But $\omega_r^* \to \omega$ as $r \uparrow \infty$, since $\int f \, d\omega_r^* = \int f \, d\omega$, for all continuous $f$ with support in $B_r$. Thus $\rho(\omega, \gamma) \leq \alpha$ follows from the lower semi-continuity of $\rho(\cdot, \gamma)$.

(iv). Let $(\omega_n)_{n \in \mathbb{N}}$ be a given sequence converging to $\omega$ in $\Gamma_X$. We have to show that $\rho_{\omega,r}(\gamma) \leq \alpha$ if there exists an $\alpha < \infty$ such that $\sup_n \rho_{\omega_n,r}(\gamma) \leq \alpha$. As above it follows that there are $\omega_n^*$ with $(\omega_n^*)_{B_r} = (\omega_n)_{B_r}$ and $\rho_{\omega_n,r}(\gamma) = \rho(\omega_n^*, \gamma) \leq \alpha$. By Lemma 4.1 $\{\omega_n^* \mid n \in \mathbb{N}\}$ is relatively compact and any accumulation point $\omega^*$ satisfies $\rho(\omega^*, \gamma) \leq \alpha$. Moreover, if $f$ is a continuous function with compact support in $B_r$,

$$\int f \, d\omega^* = \lim_{n \uparrow \infty} \int f \, d\omega_n^* = \lim_{n \uparrow \infty} \int f \, d\omega_n = \int f \, d\omega.$$

Thus $\omega^*_{B_r} = \omega_{B_r}$, and (iv) is proved. $\square$

## 5. The $\rho$-geometry on $\Gamma_X$

In this section, we will derive several auxiliary lemmas of geometric kind, which are needed in order to prove our theorems. One of the key results in this section will be Proposition 5.4 below, which roughly states that the images of some $\gamma \in \Gamma_X$ under $V_0(X)$-flows are a $\rho$-dense subset of the set of all $\omega \in \Gamma_X$ which have finite $\rho$-distance with respect to $\gamma$. This will in particular imply that any measure satisfying Assumption 1.1 has full support (cf. Proposition 5.6), a property which might be of independent interest in case of the Gibbs measures of Section 2.

Suppose we are given a $\rho$-continuous path $\xi : \mathbb{R} \to \Gamma_X$ and two real numbers $a$ and $b$ with $a < b$. Then, as usual, we define the $\rho$-energy of $\xi$ in the interval $[a, b]$ as follows.

$$E_{a,b}(\xi) = \sup\left\{ \frac{1}{2} \sum_{i=1}^n \frac{\rho(\xi_{t_i}, \xi_{t_{i-1}})^2}{t_i - t_{i-1}} \;\bigg|\; a = t_0 < t_1 < \cdots < t_n = b, \; n \in \mathbb{N} \right\}.$$

**Lemma 5.1:** *If $V \in V_0(X)$ has the flow $(\psi_t)_{t\in\mathbb{R}}$ and $\xi_t$ is given by $\xi_t = \psi_t^* \gamma$, for some fixed $\gamma \in \Gamma_X$, then the path $\xi$ is $\rho$-continuous, and*

$$E_{a,b}(\xi) = \frac{1}{2} \int_a^b \|V\|_{\xi_t}^2 \, dt,$$

*for all $a < b$.*

**Proof:** To prove '$\leq$' we remark that

$$\rho(\xi_s, \xi_t)^2 \leq \int d(\psi_s(x), \psi_t(x))^2 \gamma(dx) \leq \int \max_{s \leq r \leq t} g_{\psi_r(x)}(V, V)(t-s)^2 \, \gamma(dx).$$

Now choose an ordered partition $\Delta = \{t_0, \ldots, t_n\}$ of $[a, b]$. Then

$$\sum_{i=1}^n \frac{\rho(\xi_{t_i}, \xi_{t_{i-1}})^2}{t_i - t_{i-1}} \leq \int \sum_{i=1}^n \max_{t_{i-1} \leq r \leq t_i} g_{\psi_r(x)}(V, V) \, (t_i - t_{i-1}) \, \gamma(dx).$$

The latter converges to

$$\int \int_a^b g_{\psi_t(x)}(V, V) \, dt \, \gamma(dx) = \int_a^b \int g(V, V) \, d\xi_t \, dt$$

as the mesh of $\Delta$ tends to 0. Let us now prove that $E_{a,b}(\xi) \geq \frac{1}{2} \int_a^b \|V\|_{\xi_t}^2 \, dt$. To this end, let $E_{a,b}^X(c_x)$ denote the usual Riemannian energy of the $X$-valued curve

$c_x(t) := \psi_t(x)$, $t \in \mathbb{R}$, over the interval $[a, b]$. Then, if $\varepsilon > 0$ is given, there is some $n = n(x) \in \mathbb{N}$ such that, for $t_i := a + (b-a) \cdot i 2^{-n}$,

$$\frac{1}{2} \int_a^b g_{c_x(t)}(V, V) \, dt = E_{a,b}^X(c_x) \leq \frac{1}{2} \sum_{i=1}^{2^n} \frac{d(c_x(t_i), c_x(t_{i-1}))^2}{t_i - t_{i-1}} + \varepsilon.$$

Since $V$ has compact support, $n$ can even be chosen uniformly in $x \in \operatorname{supp} V \cap \operatorname{supp} \gamma$. Moreover, it is easy to see that, for large $n$,

$$\rho(\xi_{t_i}, \xi_{t_{i-1}}) = \left( \int d(c_x(t_i), c_x(t_{i-1}))^2 \gamma(dx) \right)^{1/2}.$$

Thus

$$E_{a,b}(\xi) \geq \frac{1}{2} \sum_{i=1}^{2^n} \frac{\rho(\xi_{t_i}, \xi_{t_{i-1}})^2}{t_i - t_{i-1}} \geq \frac{1}{2} \int_a^b \|V\|_{\xi_t} \, dt - \varepsilon \cdot \gamma(\operatorname{supp} V)$$

This proves the assertion. □

**Lemma 5.2:** *Suppose $\xi$ is as in Lemma 5.1 and $u : \Gamma_X \to \mathbb{R}$ is $\rho$-Lipschitz continuous. Then $t \mapsto u(\xi_t)$ is absolutely continuous and*

$$\left| \frac{d}{dt} u(\xi_t) \right| \leq \operatorname{Lip}(u) \cdot \|V\|_{\xi_t} \qquad \text{for almost every } t.$$

**Proof:** If $\Delta = \{t_0, t_1, \ldots, t_n\}$ is an ordered partition of some finite interval $[a, b]$, then

$$\sum_{i=1}^n \frac{(u(\xi_{t_i}) - u(\xi_{t_{i-1}}))^2}{t_i - t_{i-1}} \leq 2 \operatorname{Lip}(u)^2 E_{a,b}(\xi) = \operatorname{Lip}(u)^2 \int_a^b \|V\|_{\xi_t}^2 \, dt.$$

The lemma in Chapter II, No. 36 of Riesz and Nagy (1956) states that, under this condition, $t \mapsto u(\xi_t)$ is absolutely continuous, and that

$$\int_a^b \left( \frac{d}{dt} u(\xi_t) \right)^2 dt \leq \operatorname{Lip}(u)^2 \int_a^b \|V\|_{\xi_t}^2 \, dt.$$

Thus the lemma is proved, because $a$ and $b$ were arbitrary. □

Note that smoothness of $(X, g)$ is essential in the proof of the following lemma, and in fact the assertion may fail if $X$ is only a Lipschitz manifold.

**Lemma 5.3:** *Let $\gamma$, $\omega$, and $\eta^*$ be as in Lemma 4.1 (iv). Suppose that $(x_1, y_1)$ and $(x_2, y_2)$ are two points in the support of $\eta^*$ such that $x_1$, $y_1$, $x_2$, and $y_2$ are mutually distinct. Let $\tau_i$ denote $d(x_i, y_i)$, and assume that $c_i : [0, \tau_i] \to X$ is a minimal geodesic connecting $x_i$ with $y_i$, and suppose that $c_i$ is parameterized by arc length ($i = 1, 2$). Then neither $c_1([0, \tau_1]) \subset c_2([0, \tau_2])$ nor $c_2([0, \tau_2]) \subset c_1([0, \tau_1])$ can occur and one of the two following cases holds.*

*(i) There exists a local geodesic $c$ extending both $c_1$ and $c_2$.*

*(ii) $c_1([0, \tau_1]) \cap c_2([0, \tau_2]) = \emptyset$.*

**Proof:** Suppose that $c_1([0, \tau_1]) \subset c_2([0, \tau_2])$. But then we must have that

$$d(x_1, y_2)^2 + d(x_2, y_1)^2 < d(x_1, y_1)^2 + d(x_2, y_2)^2, \tag{5.1}$$

which contradicts the minimality of $\eta^*$.

Now suppose that $x_1$, $y_1$, $x_2$, and $y_2$ do not lie on a single local geodesic, but $c_1([0, \tau_1]) \cap c_2([0, \tau_2]) \ne \emptyset$. Then there are $t_1 \in (0, \tau_1)$ and $t_2 \in (0, \tau_2)$ such that $c_1(t_1) = c_2(t_2)$. Define two piecewise smooth curves $c_{12}$ and $c_{21}$ by

$$c_{ij}(t) = \begin{cases} c_i(t) & \text{if } t \in [0, t_i], \\ c_j(t - t_i + t_j) & \text{if } t \in [t_i, t_i + \tau_j - t_j]. \end{cases}$$

Then the total energy of $c_{12}$ and $c_{21}$ is the same as that of $c_1$ and $c_2$, i.e.,

$$E^X(c_{12}) + E^X(c_{21}) = E^X(c_1) + E^X(c_2) = d(x_1, y_1)^2 + d(x_2, y_2)^2, \tag{5.2}$$

if $E^X$ denotes the energy functional acting on piecewise smooth $X$-valued curves. Now observe that the tangent vectors of $c_1$ and $c_2$ cannot be proportional, because $c_1$ and $c_2$ cannot be extended to one single geodesic. Therefore our curves $c_{12}$ and $c_{21}$ are continuous but not differentiable in $t_1$ and $t_2$ respectively. But this implies that they cannot be energy-minimizing in the class of curves which are parameterized by arc length and have given endpoints. Hence there are two curves $c_{12}^*$ and $c_{21}^*$ parameterized by arc length connecting $x_1$ and $y_2$, and $x_2$ and $y_1$ respectively such that $E^X(c_{ij})^* < E^X(c_{ij})$. But in view of (5.2) this implies that (5.1) also holds in this case. So, again, we arrive at a contradiction to the minimality of $\eta^*$. □

**Proposition 5.4:** *Suppose $\varepsilon > 0$ is given and $\omega, \gamma \in \Gamma_X$ are such that $\rho(\omega, \gamma) < \infty$. Assume furthermore that $\gamma$ has the property that*

$$\gamma(\{x\}) \in \{0, 1\}, \qquad \text{for all } x \in X. \tag{5.3}$$

*Then, if* $\dim X \geq 2$, *there is a vector field* $V \in V_0(X)$ *with flow* $(\psi_t)_{t \in \mathbb{R}}$ *such that* $\rho(\psi_1^* \gamma, \omega) < \varepsilon$ *and* $\|V\|_{\psi_t^* \gamma} = \rho(\psi_1^* \gamma, \gamma)$, *for all* $t \in [0,1]$.

*If* $X = \mathbb{R}$, *then there are* $n \in \mathbb{N}$, $t_i \geq 0$, *and* $V_i \in V_0(X)$ *with corresponding flows* $(\psi_{i,t})_{t \in \mathbb{R}}$, $i = 1, \ldots, n$, *such that* $t_1 + \cdots + t_n = 1$ *and such that*

$$\phi_t := \psi_{i,s} \circ \psi_{i-1, t_{i-1}} \circ \cdots \circ \psi_{1, t_1}, \qquad \text{for } t = s + t_{i-1} + \cdots + t_1, \, 0 \leq s \leq t_i,$$

*satisfies* $\rho(\phi_1^* \gamma, \omega) < \varepsilon$ *and* $\rho(\phi_1^* \gamma, \gamma) \leq \|V_i\|_{\phi_t^* \gamma} \leq \rho(\phi_1^* \gamma, \gamma) + \varepsilon$, *for* $t$ *as above.*

**Proof:** First note that we can assume without loss of generality that also $\omega$ satisfies (5.3), since we can otherwise alter the corresponding points in the support of $\omega$ by an arbitrarily small portion of $\varepsilon$. Choose $r > 0$ so that the function $\rho_{\omega,r}$ of Lemma 4.2 satisfies $\rho_{\omega,r}(\gamma) \geq \rho(\omega, \gamma) - \varepsilon$. As in the proof of 4.2 (iii) there exist $\omega^*$ and $\eta^* \in \Gamma_{\gamma, \omega^*}$ such that $\omega^*_{B_r} = \omega_{B_r}$, $\rho_{\omega,r}(\gamma) = \rho(\omega^*, \gamma)$, and $\eta^*$ is optimal in the sense of 4.1 (iv) and has the form $\sum_{i=1}^{N} \delta_{(x_i, y_i)} + \sum \delta_{(x,x)}$ with $x_i \neq y_i$. Now we choose minimal geodesics $c_i : [0,1] \to X$ parameterized by arc length such that $c_i(0) = x_i$ and $c_i(1) = y_i$, $i = 1, \ldots, N$.

First consider the case $\dim X \geq 2$. Then we can assume without loss of generality that

(5.4) $$c_i([0,1]) \cap c_j([0,1]) = \emptyset \qquad \text{for } i \neq j,$$

because otherwise the situation of (i) of Lemma 5.3 must occur, and we could alter the corresponding points by a small amount to arrive at a configuration which satisfies (5.4). Now we can define $V$ to be $\dot{c}_i(t)$ in the points $c_i(t)$, for $0 \leq t \leq 1$ and $i = 1, \ldots, N$. Then the vector field $V$ is well defined due to (5.4), and we can extend it to an element in $V_0(X)$ by standard arguments. Let $(\psi_t)_{t \in \mathbb{R}}$ be its flow. Then by construction

$$\psi_t(x_i) = c_i(t), \qquad \psi_1^* \gamma = \omega^*, \qquad \text{and} \qquad \int g(V, V) \, d\psi_t^* \gamma = \rho(\gamma, \omega^*)^2,$$

for all $t \in [0,1]$ and $i = 1, \ldots, N$.

Next consider the case where $X = \mathbb{R}$. Then we are always in the situation of Lemma 5.3 (i), and we can no longer assume (5.4). However, we know that it cannot occur that $c_i([0,1]) \subset c_j([0,1])$, for some $i \neq j$, and no point $x$ such that $(x,x) \in \operatorname{supp} \eta^*$ can be contained in $\bigcup_i c_i([0,1])$. Moreover, if $c_i([0,1]) \cap c_j([0,1]) \neq \emptyset$, then $c_i$ and $c_j$ have the same orientation, i.e., $(x_i - y_i)(x_j - y_j) > 0$. Now rearrange the $x_i$s such that $(x_1, y_1), \ldots, (x_m, y_m)$ are precisely those pairs with $x_i < y_i$ and such that $x_1 < x_2 < \cdots < x_m$ and $x_{m+1} > x_{m+2} > \cdots > x_N$. Then define a piecewise constant vector field $\tilde{V}_1$ as $y_i - x_i$ on $[x_i, y_i \wedge x_{i+1})$, $i = 1, \ldots, m$, and $y_i - x_i$ on $(y_i \vee x_{i+1}, x_i]$, and $\tilde{V}_1 \equiv 0$ elsewhere. Let $(\tilde{\psi}_{1,t})_{t \in \mathbb{R}}$ denote its flow, and define

$$\tilde{t}_1 := \inf \left\{ t > 0 \,\bigg|\, \exists\, i, j \text{ such that } c_i(t) = x_j \right\} \wedge 1.$$

By construction $\tilde\psi^*_{1,t}\gamma$ is an energy minimizing $\rho$-geodesic on $[0,\tilde t_1]$ and
$$\rho(\gamma,\omega^*) = \rho(\gamma,\tilde\psi^*_{1,\tilde t_1}\gamma) + \rho(\tilde\psi^*_{1,\tilde t_1}\gamma,\omega^*).$$
In particular,
$$\eta_1^* := \int \delta_{(\tilde\psi_{1,\tilde t_1}(x),y)}\,\eta^*(dx,dy)$$
is again optimal in the sense of 4.1 (iv). Hence we can replace $\gamma$ by $\tilde\psi^*_{1,\tilde t_1}\gamma$ and $\eta^*$ by $\eta_1^*$, and start the above construction over again. Then we get a piecewise constant vector field $\tilde V_2$ with flow $(\tilde\psi_{2,t})_{t\in\mathbb R}$ and some $\tilde t_2 \in (0, 1-\tilde t_1]$ as above. It is easy to see that the above algorithm stops at some finite $n$, i.e., $\tilde t_1 + \cdots + \tilde t_n = 1$ and $\omega^* = \tilde\psi^*_{n,\tilde t_n} \circ \cdots \circ \tilde\psi^*_{1,\tilde t_1}\gamma$. Applying a mollifier to $\tilde V_1,\ldots,\tilde V_n$ gives the result. □

Let $\|\cdot\|_\infty$ denote the sup-norm of a function on $X$, and extend this norm to vector fields $V = (V^1,\ldots,V^d)$ by setting
$$\|V\|_\infty := \sup_{x\in X} \sqrt{g_x(V,V)}. \tag{5.5}$$

**Lemma 5.5:** *If $V$, $W \in V_0(X)$, $(\psi_t)_{t\in\mathbb R}$ and $(\phi_t)_{t\in\mathbb R}$ denote the corresponding flows, and $A = \operatorname{supp} V \cup \operatorname{supp} W$. Then there is a constant $c = c(V,A)$ such that*
$$\rho(\psi_t^*\gamma, \phi_t^*\gamma) \le ct \cdot e^{ct} \cdot \|V-W\|_\infty \cdot \sqrt{\gamma(A)}, \qquad \text{for all } \gamma \in \Gamma_X \text{ and } t \ge 0.$$

**Proof:** By the Nash (1956) embedding theorem $X$ can be isometrically embedded into some $\mathbb R^n$. Below $|\cdot|_1$ will denote the $\ell^1$-norm on $\mathbb R^n$. Then
$$|\psi_t(x) - \phi_t(x)|_1 \le \int_0^t \big|V(\psi_s(x)) - W(\phi_s(x))\big|_1\, ds$$
$$\le \left(t\cdot \sup_{y\in X} |V(y)-W(y)|_1 + L_V\cdot \int_0^t |\psi_s(x)-\phi_s(x)|_1\, ds\right)\cdot \mathrm{I}_A(x),$$
where, for $V = (V^1,\ldots,V^n)$,
$$L_V := \sum_{i=1}^n \sup\left\{\frac{|V^i(y)-V^i(z)|}{|z-y|_1}\ \Big|\ z\neq y,\ z,y\in X\right\}.$$
Gronwall's lemma now yields
$$|\psi_t(x) - \phi_t(x)|_1 \le t\cdot \sup_{y\in X}|V(y)-W(y)|_1 \cdot e^{tL_V}\cdot \mathrm{I}_A(x).$$
Since $A$ is compact, there are positive constants $c_1$ and $c_2$ depending only on $A$ such that, for all $y$, $z\in A$ and every tangent vector $U \in T_y X$,
$$\frac{1}{c_1}|U|_1 \le \sqrt{g_y(U,U)} \le c_1|U|_1 \qquad\text{and}\qquad d(y,z) \le c_2|y-z|_1.$$
By taking $\eta = \int \delta_{(\psi_t(x),\phi_t(x))}\,\gamma(dx)$ the assertion follows. □

**Proposition 5.6:** *Suppose that $\mu$ is a probability measure on $\Gamma_X$ which satisfies Assumptions 1.1 (a), (c), and the quasi-invariance of (d). Then $\mu$ has full support.*

**Proof:** We only give the proof in the case where $X$ is non-compact and $\dim X \geq 2$. Only minor modifications are necessary to handle the other cases. Suppose that $\mu$ does not have full support. Then we can find an open geodesic ball $B_r \subset X$, $\delta > 0$, $n \in \mathbb{N}$, and $y_1, \ldots, y_n \in B_r$ such that the open set

$$U := \left\{ \omega \in \Gamma_X \,\Big|\, \omega(\partial B_r) = 0 \text{ and } \omega_{B_r} = \sum_{i=1}^n \delta_{x_i} \text{ with } \sum_{i=1}^n d(x_i, y_i)^2 < \delta \right\}$$

satisfies $\mu(U) = 0$, because these sets form a neighborhood base for $\Gamma_X$. Now choose a compact set

$$K \subset \left\{ \gamma \in \Gamma_X \,\Big|\, \gamma(X) \geq n \text{ and } \gamma(\{x\}) \in \{0,1\} \text{ for all } x \right\}$$

with $\mu(K) > 0$. Such a set exists by Assumption 1.1 (a) and (c). Fix $\gamma \in K$. By Lemma 4.2 we must have $\rho_{\omega,r}(\gamma) < \infty$, for all $\omega \in U$. In particular, there is $\omega \in U$ with $\rho(\gamma, \omega) < \infty$. Since $U$ is open, Proposition 5.4 implies that there is even a vector field $V_\gamma \in V_0(X)$ such that its flow $(\psi_{\gamma,t})_{t \in \mathbb{R}}$ satisfies $\psi_{\gamma,1}^* \gamma \in U$. Hence $K$ is covered by $(\psi_{\gamma,1}^*)^{-1} U$, $\gamma \in K$, and we can extract a finite subcover $(\psi_{1,1}^*)^{-1} U, \ldots, (\psi_{n,1}^*)^{-1} U$. But by quasi-invariance $\mu((\psi_{i,1}^*)^{-1} U) = 0$ which contradicts $\mu(K) > 0$. □

## 6. Dirichlet forms on $\Gamma_X$.

In this section, we prove Proposition 1.4. The following two lemmas will also be needed for the proofs our main results. We suppose in this section that $\mu$ satisfies Assumption 1.1.

**Lemma 6.1:** *Let $V \in V_0(X)$ and denote its flow by $(\psi_t)_{t \in \mathbb{R}}$. Then*

$$(6.1) \qquad \int (u \circ \psi_t^* - u) v \, d\mu = \int_0^t \int (u \circ \psi_s^*) \nabla_V^{\Gamma *} v \, d\mu \, ds,$$

*for all bounded and measurable $u$ and all $v \in \mathcal{F}C_b^\infty$.*

**Proof:** By a monotone class argument it suffices to prove (6.1) for $u \in \mathcal{F}C_b^\infty$. Then also $u \circ \psi_s^* \in \mathcal{F}C_b^\infty$. Note that

$$u \circ \psi_t^* - u = \int_0^t (\nabla_V^\Gamma u) \circ \psi_s^* \, ds = \int_0^t \nabla_V^\Gamma (u \circ \psi_s^*) \, ds.$$

Therefore, by Assumption 1.1 (e),

$$\int (u \circ \psi_t^* - u) v \, d\mu = \int_0^t \int (u \circ \psi_s^*) \nabla_V^{\Gamma *} v \, d\mu \, ds.$$

This proves the lemma. □

**Lemma 6.2:** *If $u \in \mathbb{F}$ and $(\psi_t)_{t \in \mathbb{R}}$ is the flow of some $V \in V_0(X)$, then*

$$u \circ \psi_t^* - u = \int_0^t \langle \nabla^\Gamma u, V \rangle \circ \psi_s^* \, ds, \qquad \text{for all } t \in \mathbb{R}, \mu\text{-a.e.},$$

*and, for all $v \in \mathcal{F}C_b^\infty$ and all $s \in \mathbb{R}$,*

$$\int \langle \nabla^\Gamma u, V \rangle \circ \psi_s^* \cdot v \, d\mu = \int (u \circ \psi_s^*) \nabla_V^{\Gamma*} v \, d\mu. \tag{6.2}$$

**Proof:** Let $v \in \mathcal{F}C_b^\infty$. Then by definition of $\mathbb{F}$ and Lemma 6.1

$$\begin{aligned}
\int \langle \nabla^\Gamma u, V \rangle \circ \psi_s^* \cdot v \, d\mu &= \lim_{t \to 0} \int \frac{1}{t} \big( u \circ \psi_{s+t}^* - u \circ \psi_s^* \big) v \, d\mu \\
&= \lim_{t \to 0} \frac{1}{t} \int_0^t \int (u \circ \psi_{s+r}^*) \nabla_V^{\Gamma*} v \, d\mu \, dr \\
&= \int (u \circ \psi_s^*) \nabla_V^{\Gamma*} v \, d\mu,
\end{aligned}$$

because $t \mapsto u \circ \psi_t^* \in L^2(\mu)$ is continuous. Hence, again by Lemma 6.1

$$\int (u \circ \psi_t^* - u) v \, d\mu = \int \Big[ \int_0^t \langle \nabla^\Gamma u, V \rangle \circ \psi_s^* \, ds \Big] v \, d\mu.$$

By continuity of $t \mapsto u \circ \psi_t^* \in L^2(\mu)$ the assertion follows. $\square$

Following Albeverio, Kondratiev and Röckner (1997b) and Eberle (1995), we let $\mathcal{V}\mathcal{F}C_b^\infty$ denote the set of all "smooth vector fields" on $\Gamma_X$, i.e. the set of all sections $Y$ of $T\Gamma_X$ which are of the form $Y(\gamma, x) = \sum_{i=1}^n v_i(\gamma) V_i(x)$, for $n \in \mathbb{N}$, $v_i \in \mathcal{F}C_b^\infty$, and $V_i \in V_0(X)$. Then Assumption 1.1 (e) implies that, for $u \in \mathcal{F}C_b^\infty$,

$$\int \langle \nabla^\Gamma u(\gamma), Y(\gamma) \rangle_\gamma \, \mu(d\gamma) = - \int u \, \mathrm{div}_\mu^\Gamma Y \, d\mu,$$

where

$$\mathrm{div}_\mu^\Gamma Y(\gamma) := - \sum_{i=1}^n \nabla_{V_i}^{\Gamma*} v_i(\gamma).$$

Then $(\mathrm{div}_\mu^\Gamma, \mathcal{V}\mathcal{F}C_b^\infty)$ is a densely defined linear operator from $L^2(\Gamma_X \to T\Gamma_X, \mu)$ to $L^2(\mu)$, and we denote its adjoint by $(d^\mu, W^{1,2})$. Functions $u \in W^{1,2}$ are weakly differentiable in the sense that $\int \langle d^\mu u, Y \rangle \, d\mu = - \int u \, \mathrm{div}_\mu^\Gamma Y \, d\mu$ holds for all $Y \in \mathcal{V}\mathcal{F}C_b^\infty$.

**Lemma 6.3:** $\mathbb{F} \subset W^{1,2}$ and $\nabla^\Gamma = d^\mu$ on $\mathbb{F}$.

**Proof:** Let $u \in \mathbb{F}$ and $Y = \sum_{i=1}^n v_i V_i \in \mathcal{VFC}_b^\infty$. Then we get by using (6.2) for $s = 0$

$$\int \langle \nabla^\Gamma u, Y \rangle \, d\mu = \sum_{i=1}^n \int u \nabla_{V_i}^{\Gamma *} v_i \, d\mu = -\int u \operatorname{div}_\mu^\Gamma Y \, d\mu.$$

Hence the lemma follows. $\square$

**Lemma 6.4:** *Suppose $(a,b) \ni t \mapsto v(t) \in L^2(\mu)$ is differentiable at $t_0 \in (a,b)$ and let $\varphi \in C_b^1(\mathbb{R})$. Then $(a,b) \ni t \mapsto \varphi(v(t)) \in L^2(\mu)$ is differentiable at $t_0$ with derivative*

$$\left.\frac{d}{dt}\right|_{t=t_0} \varphi(v(t)) = \varphi'(v(t_0)) \left.\frac{dv(t)}{dt}\right|_{t=t_0}.$$

**Proof:** Let $\Delta(t) := (\varphi(v(t)) - \varphi(v(t_0)))/(t - t_0)$. It suffices to show that every sequence $t_n \to t_0$ has a subsequence $(t_{n_k})$ such that $\Delta(t_{n_k}) \to \varphi'(v(t_0))(dv(t)/dt)|_{t=t_0}$ in $L^2(\mu)$. But this follows easily by noting that $\{(\Delta(t_n))^2 \mid n \in \mathbb{N}\}$ is uniformly integrabel due to the Lipschitz property of $\varphi$ and by taking $(t_{n_k})$ so that $\mu$-a.e. $(v(t_{n_k}) - v(t_0))/(t_{n_k} - t_0) \to (dv(t)/dt)|_{t=t_0}$. $\square$

**Proof of Proposition 1.4:** (i). Since $(\operatorname{div}_\mu^\Gamma, \mathcal{VFC}_b^\infty)$ is densely defined, the form $\widehat{\mathcal{E}}(u,u) := \int \langle d^\mu u, d^\mu u \rangle \, d\mu$, $u \in W^{1,2}$, is closed. Therefore $(\mathcal{E}_\mu^\Gamma, \mathcal{FC}_b^\infty)$, $(\mathcal{E}_\mu^\Gamma, \mathbb{F})$, and $(\mathcal{E}_\mu^\Gamma, \mathbb{F}^{(c)})$ are all closable by Lemma 6.3 and Proposition I.3.5 of Ma and Röckner (1992). Thus $(\mathcal{E}_\mu^\Gamma, \mathcal{F}_0)$, $(\mathcal{E}_\mu^\Gamma, \mathcal{F})$, and $(\mathcal{E}_\mu^\Gamma, \mathcal{F}^{(c)})$ are symmetric closed forms, and it remains to check the contraction property to conclude that they are Dirichlet forms. But, for $u \in \mathbb{F}$, this is clear from (1.3) and Lemma 6.4, and if $u \in \mathcal{F}$ is arbitrary it then follows from Proposition II.4.10 of Ma and Röckner (1992).

Assertion (ii) is immediately implied by Theorem 1.3 if $u$ is bounded. The extension to $u \in L^2(\mu)$ then follows easily by approximating $u$ by $(-n) \vee u \wedge n$.

As for (iii), observe that the linear operators $\nabla^\Gamma$ and $d^\mu$ and hence also $\operatorname{div}_\mu^\Gamma \nabla^\Gamma$ and $\operatorname{div}_\mu^\Gamma d^\mu$ coincide on $\mathcal{FC}_b^\infty$. But if $\mu$ is a mixed Poisson measure as in the assertion, then due to Theorem 4.2 of Albeverio, Kondratiev and Röckner (1997a) and its proof we must have that $\mathcal{F} = \mathcal{F}^{(c)} = \mathcal{F}_0$.

(iv). For $u \in \mathcal{F}$, there is a sequence $(u_n)_{n \in \mathbb{N}} \subset \mathbb{F}$ converging to $u$ in $\mathcal{F}$. Thus, if $\Lambda$ denotes the $\sigma$-finite measure $\Lambda(dx, d\gamma) = \gamma(dx) \, \mu(d\gamma)$, $(\nabla^\Gamma u_n)_{n \in \mathbb{N}}$ is a Cauchy sequence in $L^2(X \times \Gamma_X \to TX, \Lambda)$. Hence there exists an element $\nabla^\Gamma u \in L^2(X \times \Gamma_X \to TX, \Lambda)$ such that $\mathcal{E}_\mu^\Gamma(u,u) = \int g_x(\nabla^\Gamma u(\gamma), \nabla^\Gamma u(\gamma)) \, \Lambda(dx, d\gamma)$. $\square$

## 7. Proof of Theorem 1.3

Let $u \in L^2(\mu)$ be $\rho$-Lipschitz continuous with $L := \text{Lip}(u)$. Fix $V \in V_0(X)$ having the flow $(\psi_t)_{t \in \mathbb{R}}$. Then consider the set

$$\Omega_V := \left\{ \gamma \in \Gamma_X \;\Big|\; \frac{u(\psi_t^*\gamma) - u(\gamma)}{t} \to G_V(\gamma) \text{ as } t \to 0, \text{ and } |G_V(\gamma)| \le L \cdot \|V\|_\gamma \right\}.$$

Then $\Omega_V$ is measurable because the existence of $\lim_{t \to 0} \frac{1}{t}\big(u(\psi_t^*\gamma) - u(\gamma)\big)$ is equivalent to the existence of the limit of $\frac{1}{r}\big(u(\psi_r^*\gamma) - u(\gamma)\big)$ for $r \to 0$, $r \in \mathbb{Q}$, because $t \mapsto u(\psi_t^*\gamma)$ is continuous by Lemma 5.2. Moreover, we claim that $\Omega_V$ has full $\mu$-measure. Indeed, we know by Lemma 5.2 that, for all $\gamma \in \Gamma_X$, the set of $s \in \mathbb{R}$ for which $\psi_s^*\gamma \in \Omega_V$ has full Lebesgue measure. Thus

$$0 = \int_0^1 ds \int \mu(d\gamma)\, \mathrm{I}_{\Omega_V^C}(\psi_s^*\gamma) = \int_0^1 ds \int_{\Omega_V^C} \mu(d\gamma)\, \frac{d\mu \circ (\psi_s^*)^{-1}}{d\mu}(\gamma).$$

But, for every $s$, $d\mu \circ (\psi_s^*)^{-1}/d\mu$ is $\mu$-a.s. strictly positive by Assumption 1.1 (d), and hence $\mu(\Omega_V^C) = 0$.

We now claim that

$$(7.1) \qquad \frac{u \circ \psi_t^* - u}{t} \longrightarrow G_V \qquad \text{in } L^2(\mu \circ (\psi_s^*)^{-1}) \text{ as } t \to 0, \text{ for all } s \in \mathbb{R}.$$

Indeed, we already know from Lemma 5.2 that

$$(7.2) \qquad \begin{aligned} \sup_{-1 \le t \le 1} \left| \frac{u(\psi_t^*\gamma) - u(\gamma)}{t} \right| &\le \sup_{-1 \le t \le 1} \frac{L}{t} \int_0^t \|V\|_{\psi_s^*\gamma}\, ds \\ &\le L \Big( \int \sup_{-1 \le s \le 1} [g(V,V) \circ \psi_s]\, d\gamma \Big)^{1/2} \\ &\le L\|V\|_\infty \sqrt{\gamma(\text{supp } V)}. \end{aligned}$$

By Assumption 1.1 (b), (7.1) thus follows from dominated convergence.

Observe that it follows from (7.2) that (6.1) holds for $u$ even though $u$ is not necessarily bounded. In view of (7.1), this implies that

$$\int G_V\, v\, d\mu = \int u\, \nabla_V^{\Gamma*} v\, d\mu, \qquad \text{for all } v \in \mathcal{F}C_b^\infty,$$

because $s \mapsto u(\psi_s^*\gamma)$ is continuous by Lemma 5.2. Then note that $V \mapsto \nabla_V^{\Gamma*} v$ is linear by Assumption 1.1 (e). Hence if $V$ can be written as $\alpha_1 V_1 + \cdots + \alpha_k V_k$ with $\alpha_i \in \mathbb{R}$ and $V_i \in V_0(X)$, $i = 1, \ldots, k$, then

$$\int G_V\, v\, d\mu = \sum_{i=1}^k \alpha_i \int u\, \nabla_{V_i}^{\Gamma*} v\, d\mu = \int \sum_{i=1}^k \alpha_i G_{V_i}\, v\, d\mu \qquad \text{for all } v \in \mathcal{F}C_b^\infty.$$

Therefore, we may conclude that

$$G_V = \sum_{i=1}^{k} \alpha_i G_{V_i} \qquad \mu\text{-a.s.} \tag{7.3}$$

Now let $\mathcal{V} \subset V_0(X)$ be a countable $\mathbb{Q}$-vector space such that, for all $V \in V_0(X)$, there exist $V_n \in \mathcal{V}$, $n \in \mathbb{N}$, such that $\|V - V_n\|_\infty \to 0$ as $n \uparrow \infty$ and all $V_n$ have compact support in a common compact subset of $X$, where the norm $\|\cdot\|_\infty$ was defined in (5.5). Such a space $\mathcal{V}$ can be easily constructed by using partitions of unity on $X$. Let $\Omega_\mathcal{V}$ denote the intersection of all sets $\Omega_V$ with $V \in \mathcal{V}$. Then $\mu(\Omega_\mathcal{V}) = 1$. Now take $\Gamma_0$ to be the set of all $\gamma \in \Omega_\mathcal{V}$ such that $V \mapsto G_V(\gamma)$ is a $\mathbb{Q}$-linear mapping on $\mathcal{V}$. By virtue of (7.3) we thus get $\mu(\Gamma_0) = 1$. Fix $\gamma \in \Gamma_0$. We have $|G_V(\gamma)| \leq L\|V\|_\gamma$, for all $V \in \mathcal{V}$. Hence by the above we can extend this mapping to a linear mapping (again denoted by $G_V(\gamma)$) defined on the whole of $V_0(X)$. Again

$$|G_V(\gamma)| \leq L\,\|V\|_\gamma \qquad \text{for all } V \in V_0(X).$$

Hence there is $\nabla^\Gamma u(\gamma) \in T_\gamma \Gamma_X$ such that $G_V(\gamma) = \langle \nabla^\Gamma u(\gamma), V \rangle_\gamma$, and $\|\nabla^\Gamma u(\gamma)\|_\gamma \leq L$. Therefore assertion (i) is proved.

The statement (ii) is already settled if $V \in \mathcal{V}$. If $V \in V_0(X)\backslash\mathcal{V}$ pick some $W$ in $\mathcal{V}$ such that $\|V - W\|_\infty \leq \varepsilon$, and let $(\phi_t)_{t\in\mathbb{R}}$ denote the flow generated by $W$. Then Lemma 5.5 yields

$$|u(\psi_t^*\gamma) - u(\phi_t^*\gamma)| \leq L\rho(\psi_t^*\gamma, \phi_t^*\gamma) \leq L\,t\,c\,\varepsilon\,e^{tc}\,\gamma(A)^{1/2},$$

where $A = \operatorname{supp} V \cup \operatorname{supp} W$ and $c$ is a constant depending only on $V$ and $A$. Therefore, if $\gamma \in \Gamma_0$,

$$\left|\frac{u(\psi_t^*\gamma) - u(\gamma)}{t} - \langle \nabla^\Gamma u(\gamma), V\rangle_\gamma\right|$$
$$\leq \varepsilon\,L\,c\,e^{tc}\,\gamma(A)^{1/2} + \|\nabla^\Gamma u(\gamma)\|_\gamma \cdot \|V - W\|_\gamma$$
$$+ \left|\frac{u(\phi_t^*\gamma) - u(\gamma)}{t} - \langle \nabla^\Gamma u(\gamma), W\rangle_\gamma\right|$$
$$\leq \varepsilon\,L\,c\,e^{tc}\,\gamma(A)^{1/2} + \varepsilon\,L\,\gamma(A)^{1/2} + o(1)$$

as $t \to 0$. This proves (ii). $\square$

## 8. Proof of Theorem 1.5

Throughout the proof, take $\mathcal{V} \subset V_0(X)$ as in the proof of Theorem 2.3, let $(\psi_{V,t})_{t\in\mathbb{R}}$ denote the flow of a vector field $V \in V_0(X)$, and suppose that $\mu$ satisfies Assumption 1.1. We will need the following simple lemma.

**Lemma 8.1:** *Suppose $(v_n)_{n \in \mathbb{N}} \subset L^2(\mu)$ converges to 0 in $L^2(\mu)$ and $V \in V_0(X)$ is a vector field. Then, for $r, t \in \mathbb{R}$, $r < t$,*

$$\int_r^t |v_n \circ \psi_{V,s}^*| \, ds \longrightarrow 0 \qquad \text{in } \mu\text{-probability as } n \uparrow \infty.$$

**Proof:** Let

$$\widetilde{\Phi}_s := \frac{d\mu \otimes ds}{d\mu \circ (\psi_{V,s}^*)^{-1} \otimes ds}, \qquad \text{and} \qquad \overline{\Phi}_s := \widetilde{\Phi}_s \circ \psi_{V,s}^*.$$

Then $(s, \gamma) \mapsto \overline{\Phi}_s(\gamma)$ is jointly measurable, and

$$\int |v_n|^2 \, d\mu = \frac{1}{t-r} \int_r^t \int |v_n \circ \psi_{V,s}^*|^2 \overline{\Phi}_s \, d\mu \, ds$$

$$\geq \frac{1}{t-r} \int \Big[ \int_r^t |v_n \circ \psi_{V,s}^*| \, ds \Big]^2 \cdot \Big[ \int_r^t \overline{\Phi}_s^{-1} \, ds \Big]^{-1} d\mu.$$

This implies the assertion, because $0 < \int_r^t \overline{\Phi}_s^{-1} \, ds < \infty$ holds $\mu$-a.s. by 1.1 (d). $\square$

Now we will show the first assertion of Theorem 1.5. To this end, suppose that $\hat{u}$ is $\rho$-continuous $\mu$-modification of a function $u \in \mathcal{F}$ with $\mathbf{\Gamma}_\mathcal{E}(u, u)(\gamma) = \|\nabla^\Gamma u(\gamma)\|_\gamma^2 \leq C^2$, for $\mu$-a.e. $\gamma \in \Gamma_X$. Choose a sequence $(u_n)_{n \in \mathbb{N}} \subset \mathbb{F}$ converging to $u$ in $\mathcal{F}$. By Lemma 6.2

$$(8.1) \qquad u_n(\psi_{V,t}^* \gamma) - u_n(\gamma) = \int_0^t \langle \nabla^\Gamma u_n, V \rangle_{\psi_{V,s}^* \gamma} \, ds, \qquad t \in \mathbb{R},$$

holds for $\mu$-a.e. $\gamma$ and all $V \in V_0(X)$ and $n \in \mathbb{N}$. Hence there is a measurable subset $\Omega_0$ of $\Gamma_X$ with $\mu(\Omega_0) = 1$ such that (8.1) holds for all $\omega \in \Omega_0$, $v \in \mathcal{V}$, $t \in \mathbb{R}$, and $n \in \mathbb{N}$. Next, by applying Lemma 8.1 with $v_n(\gamma) := \|\nabla^\Gamma u_n(\gamma) - \nabla^\Gamma u(\gamma)\|_\gamma$, a diagonalization argument implies the existence of a subsequence $(u_{n_k})_{k \in \mathbb{N}}$ and a measurable subset $\Omega_1 \subset \Omega_0$ with full $\mu$-measure such that, for all $\gamma \in \Omega_1$, $V \in \mathcal{V}$, $k \in \mathbb{N}$, $r < t$, and $s \in \mathbb{Q}$,

$$u_{n_k}(\psi_{V,s}^* \gamma) \to \hat{u}(\psi_{V,s}^* \gamma), \qquad u_{n_k}(\gamma) \longrightarrow \hat{u}(\gamma), \qquad \text{and}$$

$$\int_r^t \|\nabla^\Gamma u_{n_k} - \nabla^\Gamma u\|_{\psi_{V,s}^* \gamma} \, ds \to 0 \qquad \text{as } k \uparrow \infty.$$

Hence

$$(8.2) \qquad \hat{u}(\psi_{V,t}^* \gamma) - \hat{u}(\gamma) = \int_0^t \langle \nabla^\Gamma u, V \rangle_{\psi_{V,s}^* \gamma} \, ds$$

is true for all $\gamma \in \Omega_1$, $t \in \mathbb{Q}$, and all $V \in \mathcal{V}$. By Lemma 5.2 and the $\rho$-continuity of $\hat{u}$, the identity (8.2) extends to all $t \in \mathbb{R}$.

Now let $\Omega_2$ be the set $\Omega_1 \cap \{\gamma \in \Gamma_X \mid \gamma \text{ satisfies } (5.3)\}$. Then $\mu(\Omega_2) = 1$. If $\dim X \geq 2$, $\omega \in \Gamma_X$, and $\gamma \in \Omega_2$ are such that $\rho(\gamma, \omega) < \infty$, then Proposition 5.4 and Lemma 5.5 yield the existence of a sequence $(V_n)_{n \in \mathbb{N}} \subset \mathcal{V}$ such that

$$\rho(\psi_{V_n,1}^* \gamma, \omega) \leq \frac{1}{n} \quad \text{and} \quad \int_0^1 \|V_n\|_{\psi_{V_n,t}^* \gamma}^2 \, dt \leq \left(\rho(\gamma, \omega) + \frac{1}{n}\right)^2 \quad \text{for all } n \in \mathbb{N}.$$

From this, our assumptions, and equation (8.2) we conclude that

$$(8.3) \qquad |\hat{u}(\gamma) - \hat{u}(\omega)| \leq C \limsup_{n \uparrow \infty} \sqrt{\int_0^1 \|V_n\|_{\psi_{V_n,t}^* \gamma}^2 \, dt} = C\, \rho(\gamma, \omega),$$

with a similar reasoning if $X = \mathbb{R}$. In particular, $\hat{u}$ is $\rho$-Lipschitz continuous on $\Omega_2$ with Lipschitz constant less or equal to $C$. Next let

$$\tilde{u}(\gamma) := \sup_{\omega \in \Omega_2} \left[\hat{u}(\omega) - C\rho(\omega, \gamma)\right]$$

if the supremum is finite and $\tilde{u}(\gamma) = 0$ if not. Then $\tilde{u}$ is a $\mu$-measurable function such that $\tilde{u} = \hat{u}$ on $\Omega_2$ and such that $\mathrm{Lip}(\tilde{u}) \leq C$ (cf. McShane (1934)). This proves part (i) of our theorem.

In order to prove '$\leq$' in Theorem 1.5 (ii), fix $\omega$ and $\gamma$ as in the assertion, and consider the function $\rho_{\omega,r}$ defined in Lemma 4.2. By Lemma 4.1 (viii), $c \wedge \rho_{\omega,r}$ is $\rho$-Lipschitz continuous with $\mathrm{Lip}(c \wedge \rho_{\omega,r}) \leq 1$, for all $c, r > 0$. Hence $c \wedge \rho_{\omega,r} \in \mathbb{F}^{(c)}$ and $\mathbf{\Gamma}_\mathcal{E}(c \wedge \rho_{\omega,r}, c \wedge \rho_{\omega,r}) \leq 1$ $\mu$-a.e. by Theorem 1.3 and Lemma 4.2 (ii). But $c$ and $r$ were arbitrary, and hence Lemma 4.2 (iii) implies '$\leq$'. The inequality '$\geq$' of Theorem 1.5 (ii) follows from the next lemma.

**Lemma 8.2:** *Suppose $u \in \mathcal{F}$ has a vaguely continuous $\mu$-version $\hat{u}$ and satisfies $\mathbf{\Gamma}_\mathcal{E}(u, u) \leq 1$ $\mu$-a.s. Then $\hat{u}$ is already $\rho$-Lipschitz continuous and $\mathrm{Lip}(\hat{u}) \leq 1$*

**Proof:** Since $\hat{u}$ is $\rho$-continuous, we already know from (8.3) that, for $\mu$-a.e. $\omega \in \Gamma_X$, $\hat{u}$ has the property that

$$(8.4) \qquad |\hat{u}(\gamma) - \hat{u}(\omega)| \leq \rho(\gamma, \omega), \qquad \text{for all } \gamma \in \Gamma_X.$$

Now we will show that (8.4) holds for any given $\omega_0 \in \Gamma_X$. So let $\gamma \in \Gamma_X$ be such that $\rho(\omega_0, \gamma) < \infty$, and fix a point in $X$. Let $B_r$ denote the open geodesic ball of radius $r$ around this point, and let $\partial B_r$ denote its boundary. Then pick a sequence $0 < r_1 < r_2 < \cdots$ with $r_k \uparrow \infty$ and $\omega_0(\partial B_{r_k}) = 0$. We can write the restriction $(\omega_0)_{B_{r_k}}$ of $\omega_0$ to $B_{r_k}$ as $\sum_{i=1}^{n_k} \delta_{y^i}$. As in the proof of 5.6 we see that
(8.5)
$$U_k = \Big\{\omega \in \Gamma_X \;\Big|\; \omega(\partial B_{r_k}) = 0 \text{ and } \omega_{B_{r_k}} = \sum_{i=1}^{n_k} \delta_{z^i} \text{ with } \Big(\sum_{i=1}^{n_k} d(z^i, y^i)^2\Big)^{1/2} < 1/r_k\Big\}$$

is open. Now we choose $k_0 \in \mathbb{N}$ such that $1/r_k < \varepsilon$ and

(8.6) $$|\hat{u}(\omega) - \hat{u}(\omega_0)| < \varepsilon \qquad \text{for all } \omega \in U_{k_0}.$$

Note that (8.6) automatically holds for $k_0$ replaced by any $k \geq k_0$, because $U_k \subset U_{k_0}$. Since $\mu$ has full support by Proposition 5.6, there is an $\omega_k \in U_k$ such that $u$ satisfies (8.4) for $\omega_k$ replacing $\omega$. As in (8.5) we write $(\omega_k)_{B_{r_k}}$ as $\sum_{i=1}^{n_k} \delta_{z_k^i}$. Since $\rho(\omega_0, \gamma) < \infty$, there exists $\eta^*$ optimal in the sense of Lemma 4.1 (iv). Pick points $x^1, \ldots, x^{n_k}$ in the configuration $\gamma$ such that $(y^1, x^1), \ldots, (y^{n_k}, x^{n_k})$ are points of $\eta^*$. Defining $\gamma_k := (\omega_k)_{B_{r_k}^C} + \sum_{i=1}^{n_k} \delta_{x^i}$, we get

$$\begin{aligned}
|\hat{u}(\omega_0) - \hat{u}(\gamma_k)| &\leq |\hat{u}(\omega_0) - \hat{u}(\omega_k)| + |\hat{u}(\omega_k) - \hat{u}(\gamma_k)| \\
&\leq \varepsilon + \rho(\omega_k, \gamma_k) \\
&\leq \varepsilon + \Big(\sum_{i=1}^{n_k} d(z_k^i, x^i)^2\Big)^{1/2} \\
&\leq 2\varepsilon + \Big(\sum_{i=1}^{n_k} d(y^i, x^i)^2\Big)^{1/2} \\
&\leq 2\varepsilon + \rho(\omega_0, \gamma).
\end{aligned}$$

Finally observe that $\rho(\omega_0, \gamma) < \infty$ obviously implies that $\gamma_k \to \gamma$, vaguely as $k \uparrow \infty$. $\square$

**Acknowledgment:** A large part of this work was carried out at the Mathematical Sciences Research Institute in Berkeley, and we would like to thank the institute and the organizers of the Stochastic Analysis year, Ruth Williams and Steven N. Evans, for inviting us. We are also grateful to Byron Schmuland for providing us with the proof of Lemma 2.3 and for his careful reading of a preliminary version of the manuscript.

MICHAEL RÖCKNER, FAKULTÄT FÜR MATHEMATIK, UNIVERSITÄT BIELEFELD, D-33615 BIELEFELD, GERMANY
    roeckner@mathematik.uni-bielefeld.de

ALEXANDER SCHIED, INSTITUT FÜR MATHEMATIK, HUMBOLDT-UNIVERSITÄT, UNTER DEN LINDEN 6, D-10099 BERLIN, GERMANY
    schied@mathematik.hu-berlin.de